\documentclass[11pt]{article}%
\usepackage{amsmath}
\usepackage{amsfonts}
\usepackage{amssymb}
\usepackage{graphicx}
\usepackage{theorem}
\usepackage{makeidx}
\usepackage{verbatim}
\usepackage[colorlinks=true]{hyperref}%
\providecommand{\U}[1]{\protect\rule{.1in}{.1in}}
\newtheorem{thm}{Theorem}[section]
\newtheorem{lm}[thm]{Lemma}

\newtheorem{df}[thm]{Definition}
\newtheorem{rmk}[thm]{Remark}

{\theorembodyfont{\upshape}
\newtheorem{examp}[thm]{Example}
}
\numberwithin{equation}{section} 
\begin{document}

\title{ }

\begin{center}
\vspace*{1.5cm}

\textbf{CHAIN\ RULES\ FOR LINEAR\ OPENNESS IN\ METRIC\ SPACES. }

\textbf{APPLICATIONS\ TO\ PARAMETRIC\ VARIATIONAL\ SYSTEMS}

\vspace*{1cm}

M. DUREA

{\small {Faculty of Mathematics, "Al. I. Cuza" University,} }

{\small {Bd. Carol I, nr. 11, 700506 -- Ia\c{s}i, Romania,} }

{\small {e-mail: \texttt{durea@uaic.ro}}}

\bigskip

R. STRUGARIU

{\small Department of Mathematics, "Gh. Asachi" Technical University, }

{\small {Bd. Carol I, nr. 11, 700506 -- Ia\c{s}i, Romania,} }

{\small {e-mail: \texttt{rstrugariu@tuiasi.ro}}}
\end{center}

\bigskip

\bigskip

\noindent{\small {\textbf{Abstract: }}In this work we present a general
theorem concerning chain rules for linear openness of set-valued mappings
acting between metric spaces. As particular cases, we obtain classical and
also some new results in this field of research, including the celebrated
Lyusternik-Graves Theorem. The applications deal with the study of the
well-posedness of the solution mappings associated to parametric variational
systems. Sharp estimates for the involved regularity moduli are given.}

\bigskip


\bigskip


\bigskip

\section{Introduction}

The openness at linear rate was, over time, an intensively studied issue,
having as first landmark one of the most profound results in the theory of
linear operators -- Banach's Open Mapping Principle. Published by J. Schauder
in 1930 and by S. Banach in his famous book from 1932, it proved itself to be
very useful in a wide variety of problems, and it can be formulated in order
to emphasize the equivalent properties of linear openness and metric
regularity. The wide applicability of this principle resulted in multiple
attempts to extend it, creating another two classical results in nonlinear
analysis, for strictly differentiable functions: the Tangent Space Theorem,
proved by L.A. Lyusternik \cite{Lyu} in 1934, and the Surjection Theorem,
proved by L.M. Graves \cite{Gra} in 1950. Another landmark was the extension
of the research to the case of set-valued mappings with closed and convex
graph and this was done by Ursescu \cite{Urs1975} in 1975 and Robinson
\cite{Rob1976} in 1976, respectively. The celebrated Robinson-Ursescu Theorem
was followed by several works of Robinson and Milyutin in 1970s and 1980s
concerning the preservation of regularity (and linear openness) under
functional perturbation and the corresponding applications in the study of
generalized equations. Without being exhaustive, the following list contains
other major contributors to the further development of the field: J. P. Aubin,
A. Dontchev, H. Frankowska, A. Ioffe, A. S. Lewis, B. S. Mordukhovich, J.-P.
Penot, R. T. Rockafellar.

However, it was clearly emphasized in the works of Ioffe \cite{Ioffe2000} and
Dontchev and Rockafellar \cite{DontRock2009b} that the interrelated properties
of metric regularity, Aubin and openness at linear rate have an intimate
metric character. Generally speaking, there are two main techniques to obtain
metric regularity results. The first one, going back to the original work of
Lyusternik, is based on a constructive iterative procedure, while the second
one uses Ekeland Variational Principle. Despite the fact that Ekeland
Variational Principle works on complete metric spaces, its use in this
direction generates results stated on Banach spaces and the arguments of the
proofs are given by means of contradiction. Up to 2000s, the main literature
on this topic considers both above mentioned methods, but the general
framework is that of Banach spaces; after that, the effort to avoid the linear
structure and to work in purely metric setting became more and more desirable.

This evolution could be better observed on a concrete fundamental result - the
Lyusternik-Graves Theorem - and its successive extensions. The classical
Lyusternik-Graves Theorem establishes a metric regularity of a strictly
differentiable function from a surjectivity condition of its derivative.
Moreover, it was observed that, in fact, it can be deduced from a result of
Graves concerning the preservation of metric regularity under perturbations by
surjective continuous linear operators. Subsequently, this was extended by
Milyutin to the case where the linear perturbation is replaced by a linearly
open single-valued map. In 1996, Ursescu \cite{Urs1996} was the first to
obtain a fully set-valued extension of the above results, keeping the setting
of Banach spaces. A crucial observation is suggested by the second proof of
Theorem 6, p. 520 from the milestone survey of Ioffe (\cite{Ioffe2000}). More
precisely, it became apparent that the Lyusternik iteration process can be
successfully used when the income space is a complete metric space and the
outcome space has a linear structure with shift-invariant metric. However, the
perturbations considered by Ioffe is more restrictive than those in Ursescu's work.

Recently, following an idea of Arutyunov (\cite{Arut2007}) concerning an
extension of Nadler fixed point theorem, Ioffe \cite{Ioffe2010b} and Dontchev
and Frankowska \cite{DonFra2010} gave, on one side, openness results for
set-valued compositions, and, on the other side, fully metrical extensions,
without any linear structure, of Lyusternik-Graves Theorem.

The purpose of this work is to present a very general theorem concerning chain
rules for linear openness of set-valued mappings acting between metric spaces.
In this way we enter into a dialog with both papers \cite{Ioffe2010b} and
\cite{DonFra2010}, following the research of the authors previously developed
on Banach spaces (\cite{DurStr5}). In some particular cases, we obtain
classical and also some new results in this field of research. For instance,
the celebrated Lyusternik-Graves Theorem appears as a particular case of
composition and our approach brings into light the role of the shift-invariant
property of the metric on the outcome space. The same mechanism is shown to be
available in some situations on complete Riemannian manifolds, in relation to
another particular case of composition.

As application, we study the well-posedness of the solution mappings
associated to parametric variational systems, giving sharp estimates for the
involved regularity moduli. To this aim we introduce a local chain stability
notion which preserve the Aubin property for compositions of multifunctions.

The paper is organized as follows. After a short section of preliminaries, we
present the main result of this work which refers to the openness at linear
rate of general set-valued compositions. Then we discuss several important
particular cases where the assumptions of the main results are fulfilled.
Firstly, we take into consideration the case where the composition map is the
sum and this allows us to clearly emphasize the role of the shift-invariance
property of the metric for the formulation of the conclusion in
Lyusternik-Graves Theorem. We underline the fact that the same type of
assertions can be locally obtained without shift-invariance property but with
an appropriate change of openness rates. Secondly, we analyze the situation
when the composition map is the distance function. We identify two major cases
where this map has the desired properties: the settings of normed vector
spaces and of complete Riemannian manifolds. As a by-product we provide
another proof of a particular case of \cite[Theorem 5]{DonFra2010}.

The last section concerns some applications in the theory of parametric
variational systems. After a motivational discussion on the case of separate
variables in the definition of the composition map, we introduce a local
composition stability notion which generalizes the corresponding local sum
stability studied in \cite{DurStr5}. On this basis we investigate the metric
regularity and the Aubin property of the solution mapping associated to a very
general type of parametric variational system. We point out the boundedness
constants for the regularity moduli, extending the previous results in the
field (see \cite{ArtMord2009}, \cite{ArtMord2010}, \cite{DurStr4}).

\section{Preliminaries}

This section contains some basic definitions and results used in the sequel.
In what follows, we suppose that the involved spaces are metric spaces, unless
otherwise stated. In this setting, $B(x,r)$ and $D(x,r)$ denote the open and
the closed ball with center $x$ and radius $r,$ respectively. On a product
space we take the additive metric. If $x\in X$ and $A\subset X,$ one defines
the distance from $x$ to $A$ as $d(x,A):=\inf\{d(x,a)\mid a\in A\}.$ As usual,
we use the convention $d(x,\emptyset)=\infty.$ The excess from a set $A$ to a
set $B$ is defined as $e(A,B):=\sup\{d(a,B)\mid a\in A\}.$ For a non-empty set
$A\subset X$ we put $\operatorname*{cl}A$ for its topological closure. One
says that a set $A$ is locally complete (closed) if there exists $r>0$ such
that $A\cap D(x,r)$ is complete (closed).

Let $F:X\rightrightarrows Y$ be a multifunction. The domain and the graph of
$F$ are denoted respectively by $\operatorname*{Dom}F:=\{x\in X\mid
F(x)\neq\emptyset\}$ and $\operatorname*{Gr}F:=\{(x,y)\in X\times Y\mid y\in
F(x)\}.$ If $A\subset X$ then $F(A):=%
%TCIMACRO{\dbigcup \limits_{x\in A}}%
%BeginExpansion
{\displaystyle\bigcup\limits_{x\in A}}
%EndExpansion
F(x).$ The inverse set-valued map of $F$ is $F^{-1}:Y\rightrightarrows X$
given by $F^{-1}(y)=\{x\in X\mid y\in F(x)\}$. If $F_{1}:X\rightrightarrows
Y,F_{2}:X\rightrightarrows Z,$ we define the set-valued map $(F_{1}%
,F_{2}):X\rightrightarrows Y\times Z$ by $(F_{1},F_{2})(x):=F_{1}(x)\times
F_{2}(x).$ For a parametric multifunction $F:X\times P\rightrightarrows Y,$ we
use the notations: $F_{p}(\cdot):=F(\cdot,p)$ and $F_{x}(\cdot):=F(x,\cdot).$

\bigskip

We recall now the concepts of openness at linear rate, metric regularity and
Aubin property of a multifunction around the reference point.

\begin{df}
\label{around}Let $F:X\rightrightarrows Y$ be a multifunction and
$(\overline{x},\overline{y})\in\operatorname{Gr}F.$

(i) $F$ is said to be open at linear rate $L>0$ around $(\overline
{x},\overline{y})$ if there exist a positive number $\varepsilon>0$ and two
neighborhoods $U\in\mathcal{V}(\overline{x}),$ $V\in\mathcal{V}(\overline{y})$
such that, for every $\rho\in(0,\varepsilon)$ and every $(x,y)\in
\operatorname*{Gr}F\cap\lbrack U\times V],$%
\begin{equation}
B(y,\rho L)\subset F(B(x,\rho)). \label{Lopen}%
\end{equation}

The supremum of $L>0$ over all the combinations $(L,U,V,\varepsilon)$ for
which (\ref{Lopen}) holds is denoted by $\operatorname*{lop}F(\overline
{x},\overline{y})$ and is called the exact linear openness bound, or the exact
covering bound of $F$ around $(\overline{x},\overline{y}).$

(ii) $F$ is said to have the Aubin property (or to be Lipschitz-like) around
$(\overline{x},\overline{y})$ with constant $L>0$ if there exist two
neighborhoods $U\in\mathcal{V}(\overline{x}),$ $V\in\mathcal{V}(\overline{y})$
such that, for every $x,u\in U,$%
\begin{equation}
e(F(x)\cap V,F(u))\leq Ld(x,u). \label{LLip_like}%
\end{equation}

The infimum of $L>0$ over all the combinations $(L,U,V)$ for which
(\ref{LLip_like}) holds is denoted by $\operatorname*{lip}F(\overline
{x},\overline{y})$ and is called the exact Lipschitz bound of $F$ around
$(\overline{x},\overline{y}).$

(iii) $F$ is said to be metrically regular around $(\overline{x},\overline
{y})$ with constant $L>0$ if there exist two neighborhoods $U\in
\mathcal{V}(\overline{x}),$ $V\in\mathcal{V}(\overline{y})$ such that, for
every $(x,y)\in U\times V,$%
\begin{equation}
d(x,F^{-1}(y))\leq Ld(y,F(x)). \label{Lmreg}%
\end{equation}

The infimum of $L>0$ over all the combinations $(L,U,V)$ for which
(\ref{Lmreg}) holds is denoted by $\operatorname*{reg}F(\overline{x}%
,\overline{y})$ and is called the exact regularity bound of $F$ around
$(\overline{x},\overline{y}).$
\end{df}

The links between the previous notions are as follows (see, e.g.,
\cite[Theorem 9.43]{RocWet}, \cite[Theorems 1.52]{Mor2006}).

\begin{thm}
\label{link_around}Let $F:X\rightrightarrows Y$ be a multifunction and
$(\overline{x},\overline{y})\in\operatorname{Gr}F.$ Then $F$ is open at linear
rate around $(\overline{x},\overline{y})$ iff $F^{-1}$ has the Aubin property
around $(\overline{y},\overline{x})$ iff $F$ is metrically regular around
$(\overline{x},\overline{y})$. Moreover, in every of the previous situations,%
\[
(\operatorname*{lop}F(\overline{x},\overline{y}))^{-1}=\operatorname*{lip}%
F^{-1}(\overline{y},\overline{x})=\operatorname*{reg}F(\overline{x}%
,\overline{y}).
\]

\end{thm}

In the case of parametric set-valued maps one has the following partial
notions of linear openness, metric regularity and Aubin property around the
reference point.

\begin{df}
Let $F:X\times P\rightrightarrows Y$ be a multifunction and $((\overline
{x},\overline{p}),\overline{y})\in\operatorname{Gr}F.$

(i) $F$ is said to be open at linear rate $L>0$ with respect to $x$ uniformly
in $p$ around $((\overline{x},\overline{p}),\overline{y})$ if there exist a
positive number $\varepsilon>0$ and some neighborhoods $U\in\mathcal{V}%
(\overline{x}),$ $V\in\mathcal{V}(\overline{p}),$ $W\in\mathcal{V}%
(\overline{y})$ such that, for every $\rho\in(0,\varepsilon),$ every $p\in V$
and every $(x,y)\in\operatorname*{Gr}F_{p}\cap\lbrack U\times W],$%
\begin{equation}
B(y,\rho L)\subset F_{p}(B(x,\rho)). \label{pLopen}%
\end{equation}

The supremum of $L>0$ over all the combinations $(L,U,V,W,\varepsilon)$ for
which (\ref{pLopen}) holds is denoted by $\widehat{\operatorname*{lop}}%
_{x}F((\overline{x},\overline{p}),\overline{y})$ and is called the exact
linear openness bound, or the exact covering bound of $F$ in $x$ around
$((\overline{x},\overline{p}),\overline{y}).$

(ii) $F$ is said to have the Aubin property (or to be Lipschitz-like) with
respect to $x$ uniformly in $p$ around $((\overline{x},\overline{p}%
),\overline{y})$ with constant $L>0$ if there exist some neighborhoods
$U\in\mathcal{V}(\overline{x}),$ $V\in\mathcal{V}(\overline{p}),$
$W\in\mathcal{V}(\overline{y})$ such that, for every $x,u\in U$ and every
$p\in V,$%
\begin{equation}
e(F_{p}(x)\cap W,F_{p}(u))\leq Ld(x,u). \label{pLLip_like}%
\end{equation}

The infimum of $L>0$ over all the combinations $(L,U,V,W)$ for which
(\ref{pLLip_like}) holds is denoted by $\widehat{\operatorname*{lip}}%
_{x}F((\overline{x},\overline{p}),\overline{y})$ and is called the exact
Lipschitz bound of $F$ in $x$ around $((\overline{x},\overline{p}%
),\overline{y}).$

(iii) $F$ is said to be metrically regular with respect to $x$ uniformly in
$p$ around $((\overline{x},\overline{p}),\overline{y})$ with constant $L>0$ if
there exist some neighborhoods $U\in\mathcal{V}(\overline{x}),$ $V\in
\mathcal{V}(\overline{p}),$ $W\in\mathcal{V}(\overline{y})$ such that, for
every $(x,p,y)\in U\times V\times W,$%
\begin{equation}
d(x,F_{p}^{-1}(y))\leq Ld(y,F_{p}(x)). \label{pLmreg}%
\end{equation}

The infimum of $L>0$ over all the combinations $(L,U,V,W)$ for which
(\ref{pLmreg}) holds is denoted by $\widehat{\operatorname*{reg}}%
_{x}F((\overline{x},\overline{p}),\overline{y})$ and is called the exact
regularity bound of $F$ in $x$ around $((\overline{x},\overline{p}%
),\overline{y}).$
\end{df}

Interchanging the roles of $p$ and $x$ one gets a similar set of concepts.

\section{Linear openness of compositions}

This section is devoted to the main result of the paper, i.e. a chain rule for
linear openness of set-valued maps. In fact, the results in this section are
metric extensions of some previous assertions proved, on Banach spaces, in
\cite{DurStr5}, \cite{DurStr4}.

The starting point is an implicit multifunction theorem given (with some extra
conclusions) in \cite[Theorem 3.6]{DurStr4}. Let us remark that different
versions of this lemma are done in \cite[Theorem 3.5]{ArtMord2010} (with
functions instead of multifunctions) and in \cite[Lemma 2]{Ioffe2010b}. It
worth to be mentioned that in \cite{DurStr4} this conclusion is obtained as a
consequence of a more general implicit mapping result on Banach spaces. Here
we present the full (direct) proof, on metric spaces, for the reader's convenience.

\begin{lm}
\label{gama}Let $Y,Z,W$ be metric spaces, $G:Y\times Z\rightrightarrows W$ be
a multifunction and $({\overline{y},}\overline{z},\overline{w})\in Y\times
Z\times W$ be such that $\overline{w}\in G(\overline{y},{\overline{z}}).$
Consider next the implicit multifunction $\Gamma:Z\times W\rightrightarrows Y$
defined by%
\[
\Gamma(z,w):=\{y\in Y\mid w\in G(y,z)\}.
\]

Suppose that the following conditions are satisfied:

(i) $G$ has the Aubin property with respect to $z$ uniformly in $y$ around
$((\overline{y},{\overline{z}),}\overline{w})$ with constant $D>0;$

(ii) $G$ is open at linear rate with respect to $y$ uniformly in $z$ around
$((\overline{y},{\overline{z}),}\overline{w})$ with constant $C>0.$

Then there exists $\gamma>0$ such that, for every $(z,w),(z^{\prime}%
,w^{\prime})\in D({\overline{z},}\gamma)\times D(\overline{w},\gamma),$
\begin{equation}
e(\Gamma(z,w)\cap D(\overline{y},\gamma),\Gamma(z^{\prime},w^{\prime}%
))\leq\frac{1}{C}(Dd(z,z^{\prime})+d(w,w^{\prime})\mathbb{)}. \label{AubGam}%
\end{equation}

\end{lm}

\noindent\textbf{Proof.} The condition (i) allows us to find $\alpha>0$ such
that, for every $y\in B(\overline{y},\alpha)$ and every $z,z^{\prime}\in
B(\overline{z},\alpha),$%
\begin{equation}
e(G(y,z)\cap B(\overline{w},\alpha),G(y,z^{\prime}))\leq Dd(z,z^{\prime}).
\label{AubG}%
\end{equation}

Also, from (ii), one can choose $\alpha>0$ from before sufficiently small and
$\varepsilon>0$ such that for every $\rho\in(0,\varepsilon),$ every $z\in
B(\overline{z},\alpha)$ and every $(y,w)\in\operatorname*{Gr}G_{z}\cap\lbrack
B(\overline{y},\alpha)\times B(\overline{w},\alpha)],$%
\begin{equation}
B(w,C\rho)\subset G_{z}(B(y,\rho)). \label{opG}%
\end{equation}

Take $\gamma>0$ such that $\gamma<\dfrac{\alpha}{2},$ $\gamma<\dfrac{\alpha
}{8D},$ $\dfrac{1+\gamma}{C}(3\gamma+D\gamma)<\varepsilon,$\ and pick
arbitrary $(z,w),(z^{\prime},w^{\prime})\in D({\overline{z},}\gamma)\times
D(\overline{w},\gamma).$ Moreover, consider $y\in\Gamma(z,w)\cap
D(\overline{y},\gamma).$ Then $w\in G(y,z)\cap B(\overline{w},\dfrac{\alpha
}{2}),$ whence, by (\ref{AubG}),%
\[
d(w,G(y,z^{\prime}))\leq Dd(z,z^{\prime}).
\]

Consequently, for every $\theta\in\left(  0,\min\{\gamma,\dfrac{\alpha}%
{4}\}\right)  ,$ there exists $w^{\prime\prime}\in G(y,z^{\prime})$ such that%
\[
d(w,w^{\prime\prime})<d(w,G(y,z^{\prime}))+\theta\leq Dd(z,z^{\prime}%
)+\theta<\frac{\alpha}{2}.
\]

Therefore, $w^{\prime\prime}\in B(\overline{w},\alpha).$ If $w^{\prime\prime
}=w^{\prime},$ then $y\in\Gamma(z^{\prime},w^{\prime})$ and $d(y,\Gamma
(z^{\prime},w^{\prime}))=0.$ Suppose next that $w^{\prime\prime}%
\not =w^{\prime}.$ Then $w^{\prime}\in B(w^{\prime\prime},(1+\theta
)d(w^{\prime},w^{\prime\prime})).$ Take now $\rho_{0}:=\dfrac{1+\theta}%
{C}d(w^{\prime},w^{\prime\prime})<\dfrac{1+\gamma}{C}(3\gamma+D\gamma
)<\varepsilon$ and apply (\ref{opG}) for $z^{\prime},y,w^{\prime\prime}$ and
$\rho_{0},$ to get
\[
w^{\prime}\in B(w^{\prime\prime},(1+\theta)d(w^{\prime},w^{\prime\prime
}))\subset G(B(y,\rho_{0}),z^{\prime}).
\]

\noindent Hence, there exists $y^{\prime}\in\Gamma(z^{\prime},w^{\prime})$
with $d(y,y^{\prime})<\dfrac{1+\theta}{C}d(w^{\prime},w^{\prime\prime}).$ Then%
\begin{align*}
d(y,\Gamma(z^{\prime},w^{\prime}))  &  \leq d(y,y^{\prime})<\dfrac{1+\theta
}{C}d(w^{\prime},w^{\prime\prime})\\
&  \leq\dfrac{1+\theta}{C}(d(w^{\prime},w)+d(w,w^{\prime\prime}))\\
&  \leq\dfrac{1+\theta}{C}(d(w^{\prime},w)+Dd(z,z^{\prime})+\theta).
\end{align*}

Because $y$ was arbitrary chosen from $\Gamma(z,w)\cap D(\overline{y}%
,\gamma),$ and also, $\theta$ can be made arbitrary small, taking the supremum
for all such $y$ and making $\theta\rightarrow0$ in the above relation, one
obtains (\ref{AubGam}). $\hfill\square$

\bigskip

We are now able to formulate and to prove our main theorem, a linear openness
result for a fairly general set-valued composition. This result was proved on
Banach spaces in \cite{DurStr5}, following the usual technique offered by
Ekeland Variational Principle (see the comments in Section 1). Now, we follow
the "dual" approach, by means of iteration procedure of Lyusternik type.
Remark that we have here a purely metric result, without any linear structure,
which will allow us later to clearly emphasize the role of shift-invariance
property of the metric in the set-valued version of Lyusternik-Graves Theorem.

\begin{thm}
\label{op_comp}Let $X,Y,Z,W$ be metric spaces, $F_{1}:X\rightrightarrows Y,$
$F_{2}:X\rightrightarrows Z$ and $G:Y\times Z\rightrightarrows W$ be
multifunctions and $(\overline{x},{\overline{y}},\overline{z},\overline{w})\in
X\times Y\times Z\times W$ such that $(\overline{x},{\overline{y}}%
)\in\operatorname{Gr}F_{1},$ $(\overline{x},{\overline{z}})\in
\operatorname{Gr}F_{2}$ and $(({\overline{y},}\overline{z}),\overline{w}%
)\in\operatorname{Gr}G.$ Let $H:X\rightrightarrows W$ be given by
\[
H(x):=G(F_{1}(x),F_{2}(x))\text{ for every }x\in X.
\]
Suppose that the following assumptions are satisfied:

(i) $\operatorname{Gr}F_{1},$ $\operatorname{Gr}F_{2}$ are locally complete
around $(\overline{x},{\overline{y}}),$ $(\overline{x},\overline{z}),$
respectively, and $\operatorname{Gr}G$ is locally closed around $(({\overline
{y},}\overline{z}),\overline{w});$

(ii) $F_{1}$ is open at linear rate $L>0$ around $(\overline{x},{\overline{y}%
});$

(iii) $F_{2}$ has the Aubin property around $(\overline{x},{\overline{z}})$
with constant $M>0;$

(iv) $G$ is open at linear rate with respect to $y$ uniformly in $z$ around
$((\overline{y},{\overline{z}),}\overline{w})$ with constant $C>0;$

(v) $G$ has the Aubin property with respect to $z$ uniformly in $y$ around
$((\overline{y},{\overline{z}),}\overline{w})$ with constant $D>0;$

(vi) $LC-MD>0.$

Then there exists $\varepsilon>0$ such that, for every $\rho\in(0,\varepsilon
),$%
\[
B(\overline{w},(LC-MD)\rho)\subset H(B(\overline{x},\rho)).
\]

Moreover, there exists $\varepsilon^{\prime}>0$ such that, for every $\rho
\in(0,\varepsilon^{\prime})$ and every $(x^{\prime},y^{\prime},z^{\prime
},w^{\prime})\in B(\overline{x},\varepsilon^{\prime})\times B({\overline{y}%
,}\varepsilon^{\prime})\times{B(}\overline{z},\varepsilon^{\prime})\times
B(\overline{w},\varepsilon^{\prime})$ such that $(y^{\prime},z^{\prime}%
)\in(F_{1},F_{2})(x^{\prime})$ and $w^{\prime}\in G(y^{\prime},z^{\prime}),$
\[
B(w^{\prime},(LC-MD)\rho)\subset H(B(x^{\prime},\rho)).
\]

\end{thm}

\noindent\textbf{Proof.} Without loosing the generality, the assumptions made
upon the involved mappings yield the existence of $\alpha>0$ such that

\begin{enumerate}
\item $\operatorname{Gr}F_{1}\cap\lbrack D(\overline{x},\alpha)\times
D({\overline{y}},\alpha)]$ and $\operatorname{Gr}F_{2}\cap\lbrack
D(\overline{x},\alpha)\times D(\overline{z},\alpha)]$ are complete, and
$\operatorname{Gr}G\cap\lbrack D(\overline{y},\alpha)\times D(\overline
{z},\alpha)\times D(\overline{w},\alpha)]$ is closed.

\item for every $(x,y)\in B(\overline{x},\alpha)\times B({\overline{y}}%
,\alpha),$
\begin{equation}
d(x,F_{1}^{-1}(y))\leq\frac{1}{L}d(y,F_{1}(x)) \label{mrF1}%
\end{equation}

\item for every $x,x^{\prime}\in B({\overline{x},}\alpha),$%
\begin{equation}
e(F_{2}(x)\cap B(\overline{z},\alpha),F_{2}(x^{\prime}))\leq Md(x,x^{\prime}).
\label{lipF2}%
\end{equation}

\item for every $(z,w),(z^{\prime},w^{\prime})\in B({\overline{z},}%
\alpha)\times B(\overline{w},\alpha)$,%
\begin{equation}
e(\Gamma(z,w)\cap D(\overline{y},\alpha),\Gamma(z^{\prime},w^{\prime}%
))\leq\frac{1}{C}(Dd(z,z^{\prime})+d(w,w^{\prime})\mathbb{)}, \label{relGam}%
\end{equation}
where $\Gamma$ is defined in Lemma \ref{gama}.

\item for every $y\in B(\overline{y},\alpha)$ and every $z,z^{\prime}\in
B({\overline{z},}\alpha),$%
\begin{equation}
e(G(y,z)\cap B(\overline{w},\alpha),G(y,z^{\prime}))\leq Dd(z,z^{\prime}).
\label{lipG}%
\end{equation}

\end{enumerate}

Choose $\varepsilon>0$ such that%

\begin{equation}%
\begin{array}
[c]{l}%
\varepsilon<\alpha\\
L\varepsilon<\alpha\\
M\varepsilon<\alpha\\
(LC-MD)\varepsilon<\alpha\\
C^{-1}(1+\varepsilon)(LC-MD)\varepsilon<\alpha\\
(LC)^{-1}(1+2\varepsilon)(LC-MD)\varepsilon<\alpha\\
M(LC)^{-1}(1+3\varepsilon)(LC-MD)\varepsilon<\alpha\\
MD(LC)^{-1}(1+4\varepsilon)<1\\
\left(  1+(1+4\varepsilon)\dfrac{MD}{LC}\right)  (LC-MD)\varepsilon<\alpha
\end{array}
\label{eps}%
\end{equation}

\noindent and fix $\rho\in(0,\varepsilon).$ Take now $w\in B(\overline
{w},(LC-MD)\rho),$ $w\not =\overline{w}.$ Then, there exists $\delta
\in(0,\varepsilon)$ such that%
\begin{equation}
d(w,\overline{w})<\left(  \frac{LC}{1+4\delta}-MD\right)  \rho<\left(
\frac{LC}{1+2\delta}-MD\right)  \rho<(LC-MD)\rho. \label{delta}%
\end{equation}

Define $x_{0}:=\overline{x},$ $y_{0}:=\overline{y},$ $z_{0}:=\overline{z},$
$w_{0}:=\overline{w}.$ Using (\ref{relGam}) with $z_{0}$ instead of $z$ and
$z^{\prime},$ $w_{0}$ instead of $w$ and $w\in B(\overline{w},\alpha)$ instead
of $w^{\prime},$ one has that%

\[
d(y_{0},\Gamma(z_{0},w))\leq e(\Gamma(z_{0},w_{0})\cap B(\overline{y}%
,\alpha),\Gamma(z_{0},w))<\frac{1+\delta}{C}d(w_{0},w)
\]

\noindent(where $\delta>0$ is the one chosen before; notice that $w\not =%
w_{0}$), so there exists $y_{1}\in\Gamma(z_{0},w)$ such that%
\begin{equation}
d(y_{1},y_{0})<\frac{1+\delta}{C}d(w_{0},w)<\dfrac{1+\varepsilon}%
{C}(LC-MD)\varepsilon<\alpha. \label{y1}%
\end{equation}

Hence, $y_{1}\in B(\overline{y},\alpha).$ One can use now (\ref{mrF1}) for
$(x_{0},y_{1})$ instead of $(x,y)$ and then%
\[
d(x_{0},F_{1}^{-1}(y_{1}))\leq\frac{1}{L}d(y_{1},F_{1}(x_{0})).
\]

\noindent Because $w\not =\overline{w},$ there exists $x_{1}\in F_{1}%
^{-1}(y_{1})$ such that%
\begin{align}
d(x_{1},x_{0})  &  <d(x_{0},F_{1}^{-1}(y_{1}))+\frac{\delta}{LC}d(w_{0}%
,w)\leq\frac{1}{L}d(y_{1},F_{1}(x_{0}))+\frac{\delta}{LC}d(w_{0},w)\nonumber\\
&  \leq\frac{1}{L}d(y_{1},y_{0})+\frac{\delta}{LC}d(w_{0},w)<\frac
{1+2\varepsilon}{LC}(LC-MD)\varepsilon<\alpha. \label{x1}%
\end{align}

Consequently, $x_{1}\in B(\overline{x},\alpha).$ Now, one can use
(\ref{lipF2}) for $x_{1},x_{0}$ instead of $x,x^{\prime}$ to get that%
\[
d(z_{0},F_{2}(x_{1}))\leq e(F_{2}(x_{0})\cap B(\overline{z},\alpha
),F_{2}(x_{1}))\leq Md(x_{1},x_{0}),
\]

\noindent so there exists $z_{1}\in F_{2}(x_{1})$ such that
\begin{align}
d(z_{1},z_{0})  &  <d(z_{0},F_{2}(x_{1}))+\frac{\delta M}{LC}d(w_{0}%
,w)\nonumber\\
&  \leq Md(x_{1},x_{0})+\frac{\delta M}{LC}d(w_{0},w)<M\frac{1+3\varepsilon
}{LC}(LC-MD)\varepsilon<\alpha. \label{z1}%
\end{align}

\noindent Hence, $(y_{1},z_{1})\in(F_{1},F_{2})(x_{1})\cap\lbrack
B(\overline{y},\alpha)\times B(\overline{z},\alpha)].$

Finally, one can use (\ref{lipG}) to have that%
\[
e(G(y_{1},z_{0})\cap B(\overline{w},\alpha),G(y_{1},z_{1}))\leq Dd(z_{1}%
,z_{0}).
\]

\noindent Consequently, because $w\in G(y_{1},z_{0})\cap B(\overline{w}%
,\alpha),$ there exists $w_{1}\in G(y_{1},z_{1})$ such that%
\begin{equation}
d(w_{1},w)<d(w,G(y_{1},z_{1}))+\frac{\delta MD}{LC}d(w_{0},w)\leq
Dd(z_{1},z_{0})+\frac{\delta MD}{LC}d(w_{0},w). \label{w1}%
\end{equation}

\noindent Moreover, remark that $w_{1}\in H(x_{1}).$ Then,%
\begin{align}
d(w_{1},w)  &  <Dd(z_{1},z_{0})+\frac{\delta MD}{LC}d(w_{0},w)<MDd(x_{1}%
,x_{0})+\frac{2\delta MD}{LC}d(w_{0},w)\label{w1w}\\
&  <MD\left(  \frac{1}{L}d(y_{1},y_{0})+\frac{\delta}{LC}d(w_{0},w)\right)
+\frac{2\delta MD}{LC}d(w_{0},w)\leq\frac{MD}{LC}(1+4\delta)d(w_{0}%
,w)=Kd(w_{0},w),\nonumber
\end{align}

\noindent where $K:=\dfrac{MD}{LC}(1+4\delta)<\dfrac{MD}{LC}(1+4\varepsilon
)<1.$ Also,%
\begin{align}
d(w_{1},\overline{w})  &  \leq d(w_{1},w)+d(w,w_{0})\label{w1alfa}\\
&  \leq(1+K)(LC-MD)\rho<\left(  1+(1+4\varepsilon)\frac{MD}{LC}\right)
(LC-MD)\varepsilon<\alpha.\nonumber
\end{align}

If $w_{1}=w,$ then $w\in H(x_{1}),$ with%
\begin{align*}
d(x_{1},\overline{x})  &  \leq\frac{1}{LC}(1+2\delta)d(w_{0},w)<\frac{1}%
{LC}(1+2\delta)\left(  \frac{LC}{1+2\delta}-MD\right)  \rho\\
&  =\left(  1-(1+2\delta)\frac{MD}{LC}\right)  \rho<\rho.
\end{align*}

\noindent Hence, $w\in H(B(\overline{x},\rho))$ and the proof is finished.
Suppose next that $w_{1}\not =w.$

Now, we intend to construct the sequences $x_{n},y_{n},z_{n},w_{n}$ such that,
for $n=0,1,2,...,$ one has $w_{n}\not =w$ and%
\begin{align}
y_{n+1}  &  \in\Gamma(z_{n},w)\cap B(\overline{y},\alpha)\text{ is such that
}d(y_{n+1},y_{n})<\frac{1+\delta}{C}d(w_{n},w).\nonumber\\
x_{n+1}  &  \in F_{1}^{-1}(y_{n+1})\cap B(\overline{x},\alpha)\text{ is such
that }d(x_{n+1},x_{n})<\frac{1}{L}d(y_{n+1},y_{n})+\frac{\delta}{LC}%
d(w_{n},w)\label{sir}\\
z_{n+1}  &  \in F_{2}(x_{n+1})\cap B(\overline{z},\alpha)\text{ is such that
}d(z_{n+1},z_{n})<Md(x_{n+1},x_{n})+\frac{\delta M}{LC}d(w_{n},w),\nonumber\\
w_{n+1}  &  \in G(y_{n+1},z_{n+1})\cap B(\overline{w},\alpha)\text{ is such
that }d(w_{n+1},w)<Dd(z_{n+1},z_{n})+\frac{\delta MD}{LC}d(w_{n},w).\nonumber
\end{align}

Observe that, in view of (\ref{y1}), (\ref{x1}), (\ref{z1}), (\ref{w1}),
(\ref{w1alfa}), for $n=0,$ all the assertions from the previous formula are
satisfied by $x_{1},y_{1},z_{1},w_{1}$. Suppose that for some $p\geq1$ we have
generated $x_{1},x_{2},...,x_{p},$ $y_{1},y_{2},...,y_{p},$ $z_{1}%
,z_{2},...,z_{p}$ and $w_{1},w_{2},...,w_{p}$ satisfying (\ref{sir}).

Then%
\[
d(w_{p},w)\leq Kd(w_{p-1},w)\leq...\leq K^{p}d(w_{0},w),
\]

\noindent whence%
\begin{align}
d(w_{p},\overline{w})  &  \leq d(w_{p},w)+d(w,\overline{w})\leq(1+K^{p}%
)(LC-MD)\varepsilon\label{wpalfa}\\
&  \leq(1+K)(LC-MD)\varepsilon<(1+(1+4\varepsilon)\frac{MD}{LC}%
)(LC-MD)\varepsilon<\alpha.\nonumber
\end{align}

If $w_{p}=w,$ then $w\in H(x_{p}).$ Also, combining the inequalities from
(\ref{sir}), we have%
\begin{align}
d(x_{p},{\overline{x}})  &  \leq\overset{p}{\underset{i=1}{\sum}}%
d(x_{i},x_{i-1})\leq\frac{1+2\delta}{LC}\overset{p-1}{\underset{i=0}{\sum}%
}d(w_{i},w)\nonumber\\
&  \leq\frac{1+4\delta}{LC}\overset{p-1}{\underset{i=0}{\sum}}K^{i}%
d(w_{0},w)\leq\frac{1+4\delta}{LC}\cdot\frac{1-K^{p}}{1-K}d(w_{0}%
,w)\label{xp}\\
&  <\frac{1+4\delta}{LC\left(  1-(1+4\delta)\dfrac{DM}{LC}\right)  }\left(
\frac{LC}{1+4\delta}-DM\right)  \rho=\rho.\nonumber
\end{align}

\noindent Hence, $x_{p}\in B({\overline{x},}\rho)$ and the proof is finished.
Suppose next that $w_{p}\not =w.$

Because $w_{p},w\in B(\overline{w},\alpha),$ $y_{p}\in\Gamma(z_{p},w_{p})\cap
B(\overline{y},\alpha),$ one can find, using again (\ref{relGam}), $y_{p+1}%
\in\Gamma(z_{p},w)$ such that
\begin{equation}
d(y_{p+1},y_{p})<\frac{1+\delta}{C}d(w_{p},w). \label{yp+1}%
\end{equation}

Also,%

\begin{align}
d(y_{p+1},\overline{y})  &  \leq\overset{p+1}{\underset{i=1}{\sum}}%
d(y_{i},y_{i-1})\leq\frac{1+\delta}{C}\overset{p}{\underset{i=0}{\sum}}%
d(w_{i},w)\label{yp+1alfa}\\
&  \leq\frac{1+4\delta}{C}\overset{p}{\underset{i=0}{\sum}}K^{i}d(w_{0}%
,w)\leq\frac{1+4\delta}{C}\cdot\frac{1-K^{p+1}}{1-K}d(w_{0},w)\nonumber\\
&  \leq\frac{1+4\delta}{C(1-K)}\left(  \frac{LC}{1+4\delta}-MD\right)
\varepsilon=L\varepsilon<\alpha.\nonumber
\end{align}

\noindent Hence, $y_{p+1}\in B(\overline{y},\alpha).$ Then, one can use
(\ref{mrF1}) for $(x_{p},y_{p+1})$ instead of $(x,y)$ to have%
\[
d(x_{p},F_{1}^{-1}(y_{p+1}))\leq L^{-1}d(y_{p+1},F_{1}(x_{p})).
\]

\noindent Again, we are in the case $w\not =w_{p},$ so there exists
$x_{p+1}\in F_{1}^{-1}(y_{p+1})$ such that%
\begin{align}
d(x_{p+1},x_{p})  &  <d(x_{p},F_{1}^{-1}(y_{p+1}))+\frac{\delta}{LC}%
d(w_{p},w)\leq\frac{1}{L}d(y_{p+1},F_{1}(x_{p}))+\frac{\delta}{LC}%
d(w_{p},w)\nonumber\\
&  \leq\frac{1}{L}d(y_{p+1},y_{p})+\frac{\delta}{LC}d(w_{p},w). \label{xp+1}%
\end{align}

As above, one has%
\[
d(x_{p+1},{\overline{x}})<\rho<\varepsilon<\alpha.
\]

\noindent Hence, one can apply (\ref{lipF2}) for $x_{p+1},x_{p}$ instead of
$x,x^{\prime}$ to get $z_{p+1}\in F_{2}(x_{p+1})$ such that
\begin{equation}
d(z_{p+1},z_{p})<Md(x_{p+1},x_{p})+\frac{\delta M}{LC}d(w_{p},w). \label{zp+1}%
\end{equation}

But, using (\ref{xp}),%

\begin{align}
d(z_{p+1},\overline{z})  &  \leq\overset{p+1}{\underset{i=1}{\sum}}%
d(z_{i},z_{i-1})\leq M\overset{p+1}{\underset{i=1}{\sum}}d(x_{i}%
,x_{i-1})+\frac{\delta M}{LC}\overset{p+1}{\underset{i=1}{\sum}}%
d(w_{i-1},w)\label{zp+1alfa}\\
&  \leq\frac{(1+4\delta)M}{LC}\overset{p}{\underset{i=0}{\sum}}d(w_{i},w)\leq
M\rho<M\varepsilon<\alpha.\nonumber
\end{align}

\noindent So, $z_{p+1}\in B(\overline{z},\alpha),$ whence one can use
(\ref{lipG}) to get $w_{p+1}\in G(y_{p+1},z_{p+1})$ such that%
\begin{equation}
d(w_{p+1},w)<Dd(z_{p+1},z_{p})+\frac{\delta MD}{LC}d(w_{p},w).\label{wp+1}%
\end{equation}

\noindent Finally, one can prove like in (\ref{wpalfa})\ that $d(w_{p+1}%
,\overline{w})<\alpha$. Remark also that $w_{p+1}\in H(x_{p+1}).$

In this moment, we have completely finished the induction step, hence
(\ref{sir}) holds for every positive integer $n$.

We intend to prove next that the sequences $(x_{n}),(y_{n}),(z_{n})$ satisfy
the Cauchy condition. For this, observe first that, for every $n\in
\mathbb{N},$%
\begin{align}
d(w_{n+1},w)  &  <Dd(z_{n+1},z_{n})+\frac{\delta MD}{LC}d(w_{n},w)<MDd(x_{n+1}%
,x_{n})+\frac{2\delta MD}{LC}d(w_{n},w)\label{wp+1w}\\
&  <MD\left(  \frac{1}{L}d(y_{n+1},y_{n})+\frac{\delta}{LC}d(w_{n},w)\right)
+\frac{2\delta MD}{LC}d(w_{n},w)\leq Kd(w_{n},w)\leq K^{n}d(w_{0},w),\nonumber
\end{align}

\noindent hence $w_{n}\rightarrow w$ (because $K<1$).

Also, for every $p\in\mathbb{N}$,%
\begin{align*}
d(x_{n+p},x_{n})  &  \leq\overset{p}{\underset{i=1}{\sum}}d(x_{n+i}%
,x_{n+i-1})\leq\frac{1+2\delta}{LC}\overset{p}{\underset{i=1}{\sum}%
}d(w_{n+i-1},w)\\
&  \leq\frac{1+2\delta}{LC}\overset{p}{\underset{i=1}{\sum}}K^{n+i-2}%
d(w_{0},w)\leq\frac{1+2\varepsilon}{LC}\cdot\frac{K^{n-1}}{1-K}d(w_{0},w),
\end{align*}

\noindent so, for $n$ sufficiently large, $d(x_{n+p},x_{n})$ can be made
arbitrary small. Similar assertions hold for $(y_{n})$ and $(z_{n}),$ which
can be proven in the line of (\ref{yp+1alfa}) and (\ref{zp+1alfa}). Because
$(x_{n},y_{n},z_{n})\in\operatorname{Gr}(F_{1},F_{2})\cap\lbrack
D({\overline{x},}\alpha)\times D(\overline{y},\alpha)\times D(\overline
{z},\alpha)]$ for every $n\in\mathbb{N}$, one can find, using 1., that there
exist $(x,y,z)\in\operatorname{Gr}(F_{1},F_{2})$ such that $(x_{n},y_{n}%
,z_{n})\rightarrow(x,y,z).$ Also, because $\operatorname{Gr}G\ni(y_{n}%
,z_{n},w_{n})\rightarrow(y,z,w)$ and $\operatorname{Gr}G$ is closed, we obtain
that $w\in G(y,z),$ hence $w\in H(x).$ To complete the proof, it remains to
prove that $d(x,{\overline{x}})<\rho.$

For this, taking into account (\ref{xp}) and (\ref{delta}), observe that%
\begin{align*}
d(x,\overline{x})  &  \leq d(x,x_{n})+d(x_{n},\overline{x})\leq d(x,x_{n}%
)+\frac{1+2\delta}{LC}\cdot\frac{1-K^{n}}{1-K}d(w_{0},w)\\
&  \leq d(x,x_{n})+\frac{1+2\delta}{LC\left(  1-(1+4\delta)\dfrac{DM}%
{LC}\right)  }\left(  \frac{LC}{1+4\delta}-DM\right)  \rho\\
&  =d(x,x_{n})+\frac{1+2\delta}{1+4\delta}\rho\\
&  =d(x,x_{n})+\rho-\frac{2\delta}{1+4\delta}\rho.
\end{align*}

Since $x_{n}\rightarrow x,$ for $n$ sufficiently large, $d(x,x_{n}%
)-\dfrac{2\delta}{1+4\delta}\rho<0,$ whence $d(x,\overline{x})<\rho$ and the
proof of the first part is done.

For the second part, take $\varepsilon>0$ such that all the inequalities from
(\ref{eps}) are satisfied with $\alpha$ replaced by $\dfrac{\alpha}{2}.$
Furthermore, define $\varepsilon^{\prime}:=\dfrac{\varepsilon}{2}$ and pick
$(x^{\prime},y^{\prime},z^{\prime},w^{\prime})\in B(\overline{x}%
,\varepsilon^{\prime})\times B({\overline{y},}\varepsilon^{\prime})\times
{B(}\overline{z},\varepsilon^{\prime})\times B(\overline{w},\varepsilon
^{\prime})$ such that $(y^{\prime},z^{\prime})\in(F_{1},F_{2})(x^{\prime})$,
$w^{\prime}\in G(y^{\prime},z^{\prime}).$ Also, choose $w\in B(w^{\prime
},(LC-MD)\rho).$ Then $B(x^{\prime},\dfrac{\varepsilon}{2})\subset
B({\overline{x},}\varepsilon)\subset B({\overline{x},}\alpha),$ and similar
assertions hold for the other balls. Then the proof becomes very similar to
the one of the first part, starting with $x_{0}:=x^{\prime},y_{0}:=y^{\prime
},z_{0}:=z^{\prime},w_{0}:=w^{\prime}.$\hfill$\square$

\bigskip

We remark that the proof based on iteration procedure allows to use more
refined completeness conditions on the initial data (compare with
\cite[Theorem 3.3]{DurStr5}). The next sections are two fold: firstly, they
illustrate the main result, Theorem \ref{op_comp}, by some direct consequences
and, secondly, some applications to parametric variational systems are derived.

\section{Consequences}

We start the discussion on the consequences of Theorem \ref{op_comp} by the
study of two particular cases of set-valued maps $G$ which fulfill the
properties required in that theorem.

\begin{rmk}
\label{rmk_sum}First, we consider the situation where $Y=Z=W$ is a linear
metric space endowed with a shift-invariant metric $d$ and $G:Y\times
Y\rightrightarrows Y$ is defined by $G(y,z)=\{y-z\}.$ Take $\overline
{y},\overline{z}\in Y$ and $\overline{w}=\overline{y}-\overline{z}.$ To show
that $G$ has the Aubin property with respect to $z$ uniformly in $y$ around
$(\left(  \overline{y},\overline{z}\right)  ,\overline{w})$ with constant
$D=1$ is an easy task because for every $\varepsilon>0,$ $y\in B(\overline
{y},\varepsilon),$ $z_{1},z_{2}\in B(\overline{z},\varepsilon),$
\[
d(y-z_{1},y-z_{2})=d(z_{1},z_{2}).
\]
Similarly, in order to show that $G$ is open at linear rate with respect to
$y$ uniformly in $z$ around $(\left(  \overline{y},\overline{z}\right)
,\overline{w}),$ with constant $C=1,$ consider the next remark. Take $\rho>0$
and $w\in B(\overline{w},\rho).$ Hence%
\[
d(w,\overline{y}-\overline{z})<\rho
\]
and taking again into account the shift-invariance of $d$,%
\[
d(\overline{y},w+\overline{z})<\rho,
\]
whence $w+\overline{z}\in B(\overline{y},\rho),$ i.e. $w\in B(\overline
{y},\rho)-\overline{z}=G_{\overline{z}}(B(\overline{y},\rho)).$ The linearity
ensures the announced property.
\end{rmk}

Next, we would like to emphasize that if one drops the shift-invariance of the
metric $d,$ then the properties of $G$ cannot be guaranteed.

\begin{examp}
\label{ex_sum}Consider, for instance, $Y=\mathbb{R}$ (the real line) with the
metric $d(x,y)=\left\vert x^{3}-y^{3}\right\vert .$ This metric does not
fulfill the shift-invariance property, and $G$ is not open at linear rate
around any $(\left(  \overline{y},\overline{z}\right)  ,\overline{w}).$
Indeed, with the above notations let us take $\overline{y}=\overline{z}=1.$ If
$G$ would satisfy the openness property around $(1,1,0),$ then for any
$\rho>0$ small enough, the inequality
\[
\left\vert w^{3}\right\vert <C\rho\Leftrightarrow\left\vert w\right\vert
<\sqrt[3]{C\rho}%
\]
should imply
\begin{align*}
\left\vert (w+1)^{3}-1\right\vert  &  <\rho\Longleftrightarrow\\
\left\vert w^{3}+3w^{2}+3w\right\vert  &  =\left\vert w\right\vert
(w^{2}+3w+3)<\rho.
\end{align*}
Take $\rho\in\left(  n^{-1},2n^{-1}\right)  $ for $n\in\mathbb{N}$ large
enough and $w_{n}:=\sqrt[3]{Cn^{-1}}.$ Then one should have
\[
\left\vert w_{n}^{3}+3w_{n}^{2}+3w_{n}\right\vert <2n^{-1}%
\]
i.e.%
\[
Cn^{-1}+3\sqrt[3]{C^{2}n^{-2}}+3\sqrt[3]{Cn^{-1}}<2n^{-1}%
\]
for all $n$ large enough, which is not possible.

Nevertheless, for some metrics without shift-invariance, $G$ could fulfill the
respective properties, but the constants $C$ and $D$ would depend both on the
metric and the point $(\overline{y},\overline{z}).$ Let us consider an example
of this type. Take $Y=\mathbb{R}$ with the metric%
\[
d(x,y)=\left\{
\begin{array}
[c]{l}%
2\left\vert x-y\right\vert \text{ if }x,y\leq0\\
\left\vert x-y\right\vert \text{ if }x,y>0\\
\left\vert 2x-y\right\vert \text{ if }x\leq0,y>0\\
\left\vert x-2y\right\vert \text{ if }x>0,y\leq0.
\end{array}
\right.
\]
Obviously, this is not a shift-invariance metric. Take $\overline{y}=1,$
$\overline{z}=2.$ Then $\overline{w}=-1.$ It is easy to see that $G$ satisfy
the needed properties for $C=2$ and $D=2$ because taking small balls one works
in the left-hand sides with negative numbers and, in the right-hand sides
(after translation) with positive numbers.
\end{examp}

A second important remark of this section is as follows.

\begin{rmk}
\label{rmk_dist}Take $Y=Z$ as a metric space, $W=\mathbb{R}$ and $G:Y\times
Y\rightrightarrows\mathbb{R},$ $G(y,z)=\{d(y,z)\}.$ Take $\overline
{y},\overline{z}\in Y$ and $\overline{w}=d(\overline{y},\overline{z}).$ The
fact that $G$ has the Aubin property with respect to $z$ uniformly in $y$
around $(\left(  \overline{y},\overline{z}\right)  ,\overline{w})$ with
constant $D=1$ is immediate.

Suppose now that $\overline{w}\neq0,$ i.e. $\overline{y}\neq\overline{z}$ and
write down how the $C-$linear openness property of $G$ looks like in this
case. Denote $\alpha:=d(\overline{y},\overline{z})>0.$ Take $\varepsilon
:=2^{-1}C^{-1}\alpha,$ $\rho\in(0,\varepsilon),$ $z\in B(\overline{z}%
,4^{-1}\alpha)$, $(y,w)\in\operatorname*{Gr}G_{z}\cap\lbrack B(\overline
{y},4^{-1}\alpha)\times B(\overline{w},\alpha)]$ and $\mu\in B(w,C\rho).$ Then
$d(y,z)\geq2^{-1}\alpha>C\rho,$ because, otherwise, one gets the contradiction%
\[
d(\overline{y},\overline{z})\leq d(\overline{y},y)+d(y,z)+d(z,\overline
{z})<\alpha.
\]
One has $\left\vert \mu-d(y,z)\right\vert <C\rho,$ whence%
\[
\mu\in(d(y,z)-C\rho,d(y,z)+C\rho).
\]
But $G_{z}(B(y,\rho))=\{d(u,z)\mid u\in B(y,\rho)\}.$ Now, in order to have
the openness property for $G$ one needs the following inclusion%
\begin{equation}
(d(y,z)-C\rho,d(y,z)+C\rho)\subset\{d(u,z)\mid u\in B(y,\rho)\}.
\label{dist_surj}%
\end{equation}
Notice that the reverse inclusion always holds for $C=1$
\end{rmk}

There are some particular situations where (\ref{dist_surj}) fails. In order
to illustrate this, we provide some examples.

\begin{examp}
A very simple case is that of the discrete metric $d(x,y)=0$ if $x=y$ and
$d(x,y)=1$ otherwise. A more elaborated example is the following one. Take
$Y=\mathbb{R}^{2}$, denote by $d_{E}$ the usual (Euclidean) distance on $Y$
and consider the following distance%
\[
d(x,y)=\left\{
\begin{array}
[c]{l}%
d_{E}(x,y)\text{ if }d_{E}(x,0)=d_{E}(y,0)\\
d_{E}(x,0)+d_{E}(y,0)\text{ if }d_{E}(x,0)\neq d_{E}(y,0).
\end{array}
\right.
\]
If one takes $z=(0,0)$ and $y=(0,1),$ then $\{d(u,z)\mid u\in B(y,\rho
)\}=\{1\}$ for any $\rho<1$ since $B(y,\rho)$ is an arc of the unit circle. Of
course, in these conditions, (\ref{dist_surj}) fails.
\end{examp}

Fortunately, there are as well some general remarkable cases where the
inclusion (\ref{dist_surj}) holds.

\begin{examp}
\label{ex_dist1}Firstly, this is the case when $Y$ is a normed vector space
(endowed with a norm denoted $\left\Vert \cdot\right\Vert $)$.$ Indeed, for
any $a\in(-\rho,\rho)$ consider $u:=y+a\left\Vert y-z\right\Vert ^{-1}(y-z)\in
B(y,\rho).$ Then
\[
\left\Vert u-z\right\Vert =\left\Vert y-z\right\Vert +a
\]
which shows (\ref{dist_surj}) whence the second property of $G$ in Theorem
\ref{op_comp} holds for $C=1$.
\end{examp}

\begin{examp}
\label{ex_dist2}Secondly, we emphasize that the Riemann metric on a finite
dimensional complete connected Riemannian manifold also satisfies the property
in certain situations. In this framework we experience another restriction on
the points $\overline{y},\overline{z}$ apart from the one we have met before
(i.e. $\overline{y}\neq\overline{z}$). More specifically, in the case
$\overline{y}=\overline{z}$ it is not possible, for instance, to cover the
negative part of the interval $(d(\overline{y},\overline{z})-C\rho
,d(\overline{y},\overline{z})+C\rho)=(-C\rho,C\rho).$ In some cases, on
Riemannian manifolds, it is not possible to cover the positive part of this
interval. Roughly speaking, in order to emphasize the difficulties which could
arise, let us take the case of $S^{2}\subset\mathbb{R}^{3}$ sphere, where the
points $\overline{y},\overline{z}$ are antipodal. In this case, it is not
possible to get a distance larger than $d(\overline{y},\overline{z})=\pi.$ In
fact, it happens that for some points in $B(\overline{y},\rho),$ there are
geodesic arcs connecting these points with $\overline{z}$ which are not
minimizing. This suggests that the points $\overline{y},\overline{z}$ should
be taken close enough.

After this preliminaries, let us recall the main technical fact which allows
us to formulate our conclusions. To this end, we make use of several notations
and results from \cite{DoC}.

Let $Y$ be a finite dimensional complete connected Riemannian manifold. Recall
that the Riemann metric is given by a $2-$tensor field that is symmetric and
positively definite. A Riemann metric determines an inner product and a norm
on each tangent space $T_{y}Y,$ usually written as $\left\Vert \cdot
\right\Vert _{y},$ where the subscript $y$ may be omitted. For a smooth curve
$\gamma:I\subset\mathbb{R\rightarrow}Y,$ where $I$ is an interval,
$\gamma^{\prime}(t)\in T_{\gamma(t)}Y$ for all $t\in I$ and $\left\Vert
\gamma^{\prime}(\cdot)\right\Vert \in C^{\infty}(I)$; then one can define the
length of $\gamma$ on an interval $[a,b]\subset I$ as%
\[
l(\gamma)_{a}^{b}:=\int_{a}^{b}\left\Vert \gamma^{\prime}(t)\right\Vert dt.
\]
Given two points $y,z\in Y$ one defines the Riemannian distance from $y$ to
$z$ by,
\[
d(y,z)=\inf_{\gamma}l(\gamma)_{0}^{1},
\]
where the infimum is taken over all piecewise smooth curves $\gamma
:[0,1]\rightarrow Y$ connecting $y$ and $z.$ Thus $d$ is a distance and the
topology induced by $d$ coincides with the topology of $Y$ as a differentiable
manifold. Furthermore, the celebrated Hopf-Rinow Theorem (see \cite[Theorem
2.8, p. 146]{DoC}) asserts that $(Y,d)$ is a complete metric space. Moreover,
for every $y\in Y$ the exponential map $\exp_{p}\ $is defined on all of
$T_{p}Y.$ Recall that $\exp_{p}(0)=p$ and for $v\in T_{p}Y\setminus\{0\},$
$\exp_{p}(v)$ is a point in $Y$ obtained by going out the length $\left\Vert
v\right\Vert $ starting from $p$ along a geodesic which passes through $y$
with velocity equal to $\left\Vert v\right\Vert ^{-1}v$. A geodesic $\gamma$
connecting $y$ and $z$ is called a minimizing geodesic if its arc-length is
equal to the Riemannian distance between $y$ and $z.$ A curve $\gamma
:[0,1]\rightarrow Y$ is a minimizing geodesic connecting $y$ and $z$ (i.e.
$\gamma(0)=y$ and $\gamma(1)=z$) if and only if there is $v\in T_{y}Y$ such
that $\left\Vert v\right\Vert =d(y,z)$ and $\gamma(t)=\exp_{p}(tv)$ for each
$t\in\lbrack0,1].$ In view of the equivalence between the topology of the
Riemann distance and the original topology of $Y,$ for every $y\in Y$ and for
every $\delta$ small enough, the ball $B(y,\delta)$ with respect to the
Riemann distance coincides normal ball (or geodesic ball) of center $y$ and
radius $\delta$ (see \cite[p. 70, 146]{DoC}).

Let $\overline{y}\in Y.$ According to \cite[Proposition 4.2, p. 76]{DoC} there
exists $\beta>0$ such that the geodesic ball, (denoted $B$) of center
$\overline{y}$ and radius $\beta$ is strongly convex, i.e. any two points in
that ball can be joined by a minimizing geodesic.

Now we are ready to pass to the proof of the openness property of the Riemann
distance by proving (\ref{dist_surj}) with $C=1,\ $keeping unchanged the
notations of Remark \ref{rmk_dist} and taking $\overline{z}\in B$ different
from $\overline{y}$. However, we need first to adjust (eventually) the
constant $\varepsilon$ and the neighborhood of $\overline{y}.$ For this, one
applies \cite[Theorem 3.7, p. 72]{DoC} at $\overline{y}.$ Then there exists
$\theta>0$ such that $B(\overline{y},\theta)\subset B$ is a normal
neighborhood for each of its points. Take $y\in B(\overline{y},\min
\{2^{-1}\theta,4^{-1}\alpha\}).$ Then $B(y,2^{-1}\theta)\subset B(\overline
{y},\theta)$ and $B(y,2^{-1}\theta)$ is a normal ball at $y.$ Consider now
$\varepsilon^{\prime}=\min\{\varepsilon,2^{-1}\theta\},$ $\rho\in
(0,\varepsilon^{\prime}).$ It is enough to show that for every $a\in(0,\rho),$
there exist $u_{1},u_{2}\in B(y,\rho)$ such that $d(u_{1},z)=d(y,z)-a,$ and
$d(u_{2},z)=d(y,z)+a.$ Now $B(y,a)$ is a normal ball for $y.$ We follow the
main idea of proof of Hopf-Rinow Theorem from \cite[page 147]{DoC}. Let
$x_{0}$ the point on the (compact) boundary of $B(y,a)$ where the continuous
function $d(z,\cdot)$ attaints its minimum. Then there exists $v\in T_{y}Y$
with $\left\Vert v\right\Vert =1$ such that $x_{0}=\exp_{y}(av).$ Inspecting
the proof in \cite[page 147]{DoC} one observes that it is practically shown
that for any $s\in\lbrack0,d(y,z)],$
\[
d(\exp_{y}(sv),z)=d(y,z)-s,
\]
i.e. $y,x_{0},z$ are on a geodesic. In particular, $d(\exp_{y}%
(av),z)=d(y,z)-a.$ Then we have found $u_{1}=\exp_{y}(av)\in B(y,\rho)$ having
the desired property. Now consider $u_{2}=\exp_{y}(-av)\in B(y,\rho).$ Then
since we work on a minimizing geodesic (the points are in $B$), $d(u_{2}%
,z)=d(y,z)+d(y,u_{2})=d(y,z)+a$ and the thesis is finally proved.
\end{examp}

As a first application we deduce a well-known set-valued metric version of
Lyusternik-Graves Theorem. In view of Theorem \ref{op_comp} and Remark
\ref{rmk_sum}, the proof is straightforward. Moreover, the necessity of the
shift-invariance metric is now clearly emphasized and, taking into account the
last remark in Example \ref{ex_sum}, one can even drop this requirement in
order to get local metric versions of this result with different constants.

\begin{thm}
\label{main_const}Let $X$ be a metric space, $Y$ be a linear metric space with
shift-invariant metric, $F_{1}:X\rightrightarrows Y$ and $F_{2}%
:X\rightrightarrows Y$ be multifunctions and $(\overline{x},\overline{y}%
_{1},\overline{y}_{2})\in X\times Y\times Y$ such that $(\overline
{x},\overline{y}_{1})\in\operatorname{Gr}F_{1}$ and $(\overline{x}%
,\overline{y}_{2})\in\operatorname{Gr}F_{2}.$ Suppose that the following
assumptions are satisfied:

(i) $\operatorname{Gr}F_{1}$ and $\operatorname{Gr}F_{2}$ are locally complete
around $(\overline{x},\overline{y}_{1})$ and $(\overline{x},\overline{y}%
_{2}),$ respectively;

(ii) $F_{1}$ is metrically regular around $(\overline{x},\overline{y}_{1})$
with constant $l>0;$

(iv) $F_{2}$ has the Aubin property $(\overline{x},\overline{y}_{2})$ with
constant $m>0;$

(v) $lm<1.$

Then there exists $\varepsilon>0$ such that, for every $\rho\in(0,\varepsilon
),$
\[
B(\overline{y}_{1}-\overline{y}_{2},(l^{-1}-m)\rho)\subset(F_{1}%
-F_{2})(B(\overline{x},\rho)).
\]

Moreover, there exists $\varepsilon^{\prime}>0$ such that, for every $\rho
\in(0,\varepsilon^{\prime}),$ and every $(x^{\prime},y_{1}^{\prime}%
,y_{2}^{\prime})\in\operatorname{Gr}(F_{1},F_{2})\cap\lbrack B(\overline
{x},\varepsilon)\times B({\overline{y}}_{1}{,}\varepsilon)\times
{B(\overline{y}}_{2},\varepsilon)],$%
\[
B(y_{1}^{\prime}-y_{2}^{\prime},(l^{-1}-m)\rho)\subset(F_{1}-F_{2}%
)(B(x^{\prime},\rho)).
\]

\end{thm}

\noindent\textbf{Proof.} Apply Theorem \ref{op_comp}, for the special case
where $Y=Z=W,$ $G(y,z):=\{y-z\}$ and take into account the discussion in
Remark \ref{rmk_sum}.\hfill$\square$

\bigskip

On the basis of Theorem \ref{op_comp} and Example \ref{ex_dist1} one has the
following result.

\begin{thm}
\label{LG_new}Let $X$ be a metric space and $Y$ be a normed vector space,
$F_{1}:X\rightrightarrows Y$ and $F_{2}:X\rightrightarrows Y$ be
multifunctions and $(\overline{x},\overline{y}_{1},\overline{y}_{2})\in
X\times Y\times Y$ such that $(\overline{x},\overline{y}_{1})\in
\operatorname{Gr}F_{1}$ and $(\overline{x},\overline{y}_{2})\in
\operatorname{Gr}F_{2}$ and $\overline{y}_{1}\neq\overline{y}_{2}.$ Suppose
that the following assumptions are satisfied:

(i) $\operatorname{Gr}F_{1}$ and $\operatorname{Gr}F_{2}$ are locally complete
around $(\overline{x},\overline{y}_{1})$ and $(\overline{x},\overline{y}%
_{2}),$ respectively;

(ii) $F_{1}$ is metrically regular around $(\overline{x},\overline{y}_{1})$
with constant $l>0;$

(iii) $F_{2}$ has the Aubin property $(\overline{x},\overline{y}_{2})$ with
constant $m>0;$

(iv) $lm<1.$

Then there exists $\varepsilon>0$ such that, for every $\rho\in(0,\varepsilon
),$%
\[
B(\left\Vert \overline{y}_{1}-\overline{y}_{2}\right\Vert ,(l^{-1}%
-m)\rho)\subset\{\left\Vert y_{1}-y_{2}\right\Vert \mid y_{1}\in
F_{1}(x),y_{2}\in F_{2}(x),x\in B(\overline{x},\rho)\}.
\]

Moreover, there exists $\varepsilon^{\prime}>0$ such that, for every $\rho
\in(0,\varepsilon^{\prime}),$ and every $(x^{\prime},y_{1}^{\prime}%
,y_{2}^{\prime})\in\operatorname{Gr}(F_{1},F_{2})\cap\lbrack B(\overline
{x},\varepsilon^{\prime})\times B({\overline{y}}_{1}{,}\varepsilon^{\prime
})\times{B(\overline{y}}_{2},\varepsilon^{\prime})],$%
\[
B(\left\Vert y_{1}^{\prime}-y_{2}^{\prime}\right\Vert ,(l^{-1}-m)\rho
)\subset\{\left\Vert y_{1}-y_{2}\right\Vert \mid y_{1}\in F_{1}(x),y_{2}\in
F_{2}(x),x\in B(x^{\prime},\rho)\}.
\]

\end{thm}

\bigskip

In the same line, taking into account Example \ref{ex_dist2} instead of
Example \ref{ex_dist1}, one can formulate a similar result on Riemannian manifolds.

\bigskip

Our intention is to obtain next some fixed point assertions, which are
equivalent, on Banach spaces, to the set-valued Lyusternik-Graves Theorem. In
\cite[Theorem 5]{DonFra2010}, a more general variant of the next result
concerning fixed points, stated on metric spaces, is shown to take place. In
the same paper, it is proved that this result can be used to obtain the
set-valued Lyusternik-Graves Theorem on metric spaces. The converse is also
true on Banach spaces, as proved in \cite[Theorem 4.4]{DurStr5}.

We provide next another proof of this result, when the output space is a
normed vector space.

\begin{thm}
\label{fixp}Let $X$ be a metric space, $Y$ be a normed vector space,
$F_{1}:X\rightrightarrows Y,$ $F_{2}:X\rightrightarrows Y$ be multifunctions
and $(\overline{x},\overline{y})\in X\times Y$ such that $(\overline
{x},\overline{y})\in\operatorname{Gr}F_{1}\cap\operatorname{Gr}F_{2}.$ Suppose
that the following assumptions are satisfied:

(i) $\operatorname{Gr}F_{1}$ and $\operatorname{Gr}F_{2}$ are locally complete
around $(\overline{x},\overline{y})$;

(ii) $F_{1}$ is metrically regular around $(\overline{x},\overline{y}_{1})$
with constant $l>0;$

(iii) $F_{2}$ has the Aubin property $(\overline{x},\overline{y}_{2})$ with
constant $m>0;$

(iv) $lm<1.$

Then there exist $\alpha,\beta>0$ such that, for every $x\in B(\overline
{x},\alpha),$ one has%
\begin{equation}
d(x,\operatorname*{Fix}(F_{1}^{-1}F_{2}))\leq(l^{-1}-m)^{-1}d(F_{1}(x)\cap
B(\overline{y},\beta),F_{2}(x)), \label{diffix}%
\end{equation}

\noindent where
\[
\operatorname*{Fix}(F_{1}^{-1}F_{2}):=\{x\in X\mid F_{1}(x)\cap F_{2}%
(x)\neq\emptyset\}.
\]

\end{thm}

\noindent\textbf{Proof.} Our intention is to apply Theorem \ref{LG_new} in
order to get (\ref{diffix}). Using the assumptions made, one can find
$\gamma>0$ such that the assumptions 1.-5. from the beginning of the proof of
Theorem \ref{op_comp} are satisfied with $\alpha$ replaced by $\gamma,$
$L^{-1}$ replaced by $l,$ $M$ replaced by $m,$ $C=D=1$ and $G$ as in Remark
\ref{rmk_dist}. Next, choose $\varepsilon>0$ as in (\ref{eps}), but with
$2^{-1}\gamma$ instead of $\alpha$.

Now, take $\rho\in(0,\min\{\varepsilon,(6^{-1}(l^{-1}-m)^{-1}\gamma\}),$
$\beta:=2^{-1}(l^{-1}-m)\rho$ and $\alpha>0$ such that $\alpha<\rho
,m\alpha<\beta$. Finally, fix arbitrary $x\in B(\overline{x},\alpha).$ Because
$\alpha<\rho<\varepsilon<2^{-1}\gamma$, one has that $B(x,\frac{\gamma}%
{2})\subset B(\overline{x},\gamma).$ Remark that if $F_{1}(x)\cap
B(\overline{y},\beta)=\emptyset$ or $F_{2}(x)=\emptyset$ or $F_{1}(x)\cap
B(\overline{y},\beta)\cap F_{2}(x)\not =\emptyset,$ relation (\ref{diffix})
trivially holds. Suppose next that none of the previous relations is
satisfied, and take arbitrary $y_{1}\in F_{1}(x)\cap B(\overline{y},\beta).$
Then $D(y_{1},\frac{\gamma}{2})\subset D(\overline{y},\gamma).$ Also, for
every $\mu>0,$ there exists $y_{2}^{\mu}\in F_{2}(x)$ (hence $y_{2}^{\mu
}\not =y_{1}$) such that%
\begin{equation}
0<\left\Vert y_{1}-y_{2}^{\mu}\right\Vert <d(y_{1},F_{2}(x))+\mu. \label{rel}%
\end{equation}

Using (iii), one has that%
\[
d(\overline{y},F_{2}(x))\leq e(F_{2}(\overline{x})\cap V,F_{2}(x))\leq
md(x,\overline{x}).
\]

\noindent Hence, for arbitrary $\mu>0$, there exists $v_{\mu}\in F_{2}(x)$
such that
\[
\left\Vert v-\overline{y}\right\Vert \leq md(x,\overline{x})+\mu\leq
m\alpha+\mu.
\]

Then, for $\mu>0$ sufficiently small, $\left\Vert v-{\overline{y}}\right\Vert
<\beta.$ Now, because $d(y_{1},F_{2}(x))\leq\left\Vert y_{1}-v\right\Vert
<2\beta,$ for $\mu>0$ sufficiently small, we have also that $d(y_{1}%
,F_{2}(x))+\mu<2\beta.$ Then $\left\Vert y_{1}-y_{2}^{\mu}\right\Vert <2\beta$
and, furthermore,%
\[
\left\Vert y_{2}^{\mu}-\overline{y}\right\Vert \leq\left\Vert y_{2}^{\mu
}-y_{1}\right\Vert +\left\Vert y_{1}-\overline{y}\right\Vert \leq3\beta
<2^{-1}\gamma.
\]

\noindent Consequently, for $\mu$ small enough, $D(y_{2}^{\mu},\frac{\gamma
}{2})\subset D(\overline{y},\gamma).$

Denote $\rho_{0}:=(l^{-1}-m)^{-1}[d(y_{1},F_{2}(x))+\mu]>0$ and remark that
$\rho_{0}<2(l^{-1}-m)^{-1}\beta=\rho<\varepsilon$ for $\mu$ sufficiently
small. Then all the assumptions of Theorem \ref{LG_new} are satisfied for the
reference points $x,y_{1},y_{2}^{\mu},$ so one can infer that%
\[
0\in B(\left\Vert y_{1}-y_{2}^{\mu}\right\Vert ,d(y_{1},F_{2}(x))+\mu
)\subset\{\left\Vert y_{1}^{\prime}-y_{2}^{\prime}\right\Vert \mid
y_{1}^{\prime}\in F_{1}(u),y_{2}^{\prime}\in F_{2}(u),u\in B(x,\rho_{0})\}.
\]

Hence, there exists $u\in B(x,\rho_{0})$ and $y_{1}^{\prime}\in F_{1}%
(u),y_{2}^{\prime}\in F_{2}(u)$ such that $0=\left\Vert y_{1}^{\prime}%
-y_{2}^{\prime}\right\Vert .$ Then $u\in\operatorname*{Fix}(F_{1}^{-1}F_{2})$
and%
\[
d(x,\operatorname*{Fix}(F_{1}^{-1}F_{2}))\leq d(x,u)<(l^{-1}-m)^{-1}%
[d(y_{1},F_{2}(x))+\mu].
\]

Making $\mu\rightarrow0,$ one gets (\ref{diffix}).\hfill$\square$

\bigskip

This theorem can be putted into relation with some very recent results in
literature (see \cite[Theorem 5]{DonFra2010}, \cite[Theorems 2, 3]%
{Ioffe2010b}, \cite[Theorem 4.4]{DurStr5} and the references therein).

\section{Applications}

\subsection{Local stability of compositions}

In this section we further investigate some general consequences of Theorem
\ref{op_comp} on parametric systems. To this end we recall first a stability
notion given in \cite{DurStr4}.

\begin{df}
\label{sum-sta}Let $F:X\rightrightarrows Y,$ $G:X\rightrightarrows Y$ be
multifunctions and $(\overline{x},{\overline{y},}\overline{z})\in X\times
Y\times Y$ such that ${\overline{y}\in F(\overline{x}),}$ ${\overline{z}\in
G(\overline{x}).}$ We say that the multifunction $(F,G)$ is locally sum-stable
around $(\overline{x},{\overline{y},}\overline{z})$ if for every
$\varepsilon>0$ there exists $\delta>0$ such that, for every $x\in
B(\overline{x},\delta)$ and every $w\in(F+G)(x)\cap B({\overline{y}+}%
\overline{z},\delta),$ there exist $y\in F(x)\cap B({\overline{y},}%
\varepsilon)$ and $z\in G(x)\cap B(\overline{z},\varepsilon)$ such that
$w=y+z.$
\end{df}

This definition has a parametric variant as well.

\begin{df}
\label{sum-sta_param}Let $F:X\times P\rightrightarrows Y,$
$G:X\rightrightarrows Y$ be multifunctions and $(\overline{x},\overline
{p},{\overline{y},}\overline{z})\in X\times P\times Y\times Y$ such that
${\overline{y}\in F(\overline{x},\overline{p}),}$ ${\overline{z}\in
G(\overline{x}).}$ We say that the multifunction $(F,G)$ is locally sum-stable
around $(\overline{x},\overline{p},{\overline{y},}\overline{z})$ if for every
$\varepsilon>0$ there exists $\delta>0$ such that, for every $(x,p)\in
B(\overline{x},\delta)\times B(\overline{p},\delta)$ and every $w\in
(F_{p}+G)(x)\cap B({\overline{y}+}\overline{z},\delta),$ there exist $y\in
F_{p}(x)\cap B({\overline{y},}\varepsilon)$ and $z\in G(x)\cap B(\overline
{z},\varepsilon)$ such that $w=y+z.$
\end{df}

We have now the tools in order to investigate some particular cases of Theorem
\ref{op_comp}. Consider the case where $W$ is a linear space with a
shift-invariant metric. First, take $G:Y\times Z\rightrightarrows W$ with
separate variables, i.e.
\[
G(y,z)=R(y)+T(z),
\]
where $R:Y\rightrightarrows W$ and $T:Z\rightrightarrows W.$ Observe now that
$H$ in Theorem \ref{op_comp} has the form%
\[
H(x)=R\circ F_{1}(x)+T\circ F_{2}(x).
\]
For this situation we formulate the following result.

\begin{thm}
\label{t_sep_1}Let $X,Y,Z$ be metric spaces, $W$ be a metric space with linear
structure such that the metric is shift-invariant. Take $F_{1}%
:X\rightrightarrows Y,$ $F_{2}:X\rightrightarrows Z,$ $R:Y\rightrightarrows
W,$ $T:Z\rightrightarrows W$ and $(\overline{x},\overline{y},\overline
{z},\overline{w}_{1},\overline{w}_{2})\in X\times Y\times Z\times W\times W$
such that $\overline{y}\in F_{1}(\overline{x}),$ $\overline{w}_{1}\in
R(\overline{y}),$ $\overline{z}\in F_{2}(\overline{x}),$ $\overline{w}_{2}\in
T(\overline{z})$. Suppose that:

(i) $\operatorname{Gr}F_{1},$ $\operatorname{Gr}F_{2}$ are locally complete
around $(\overline{x},{\overline{y}}),$ $(\overline{x},\overline{z}),$
respectively, and $\operatorname{Gr}G$ is locally closed around $(({\overline
{y},}\overline{z}),\overline{w});$

(ii) $F_{1}$ is open at linear rate $L>0$ around $(\overline{x},\overline
{y});$

(iii) $R$ is open at linear rate $C>0$ around $(\overline{y},\overline{w}%
_{1});$

(iv) $F_{2}$ has the Aubin property around $(\overline{x},\overline{z})$ with
constant $M>0;$

(v) $T$ has the Aubin property around $(\overline{z},\overline{w}_{2})$ with
constant $D>0;$

(vi) $(R,T)$ is locally sum-stable around $(\overline{z},{\overline{y}%
,}\overline{w}_{1},\overline{w}_{2})$ in the sense of Definition
\ref{sum-sta_param} (where consider $P=Y$ and $R$ is formally taken as
$R(y,z)=R(y),$ for all $(y,z)\in Y\times Z,$ i.e. constant with respect to $z$).

(vii) $LC-MD>0.$

Then there exists $\varepsilon>0$ such that for every $\rho\in(0,\varepsilon)$
such that%
\[
B(\overline{w}_{1}+\overline{w}_{2},(LC-MD)\rho)\subset(R\circ F_{1}+T\circ
F_{2})(B(\overline{x},\rho)).
\]

Moreover, there exists $\varepsilon^{\prime}>0$ such that for every $\rho
\in(0,\varepsilon^{\prime})$ and every $(x,y,z,w_{1},w_{2})\in B(\overline
{x},\varepsilon^{\prime})\times B(\overline{y},\varepsilon^{\prime})\times
B(\overline{z},\varepsilon^{\prime})\times B(\overline{w}_{1},\varepsilon
^{\prime})\times B(\overline{w}_{2},\varepsilon^{\prime})$ such that $y\in
F_{1}(x),$ $w_{1}\in R(y),$ $z\in F_{2}(x),$ $w_{2}\in T(z)$%
\[
B(w_{1}+w_{2},(LC-MD)\rho)\subset(R\circ F_{1}+T\circ F_{2})(B(x,\rho)).
\]

\end{thm}

\noindent\textbf{Proof.} It is enough to prove that conditions (iii), (v) and
(vi) ensure the properties of $G.$ There exist $\nu>0,$ $U\in\mathcal{V}%
(\overline{y}),$ $W_{1}\in\mathcal{V}(\overline{w}_{1})$ such that for every
$\theta\in(0,\nu)$ and every $(y,w_{1})\in\operatorname*{Gr}R\cap\lbrack
U\times W_{1}],$%
\[
B(w_{1},C\theta)\subset R(B(y,\theta)).
\]
Take $\varepsilon>0$ such that $B(\overline{w}_{1},\varepsilon)\subset W_{1}.$
Then, from (vi), there exists $\delta>0$ such that for every $(z,y)\in
B(\overline{z},\delta)\times B(\overline{y},\delta)$ and every $w\in
(R_{y}+T)(z)\cap B(\overline{w},\delta),$ there exist $w_{1}\in R(y)\cap
B({\overline{w}}_{1}{,}\varepsilon)$ and $z\in T(z)\cap B(\overline{w}%
_{2},\varepsilon)$ such that $w=w_{1}+w_{2}.$

Let $\mu>0$ such that $B(\overline{y},\mu)\subset U$ and consider $U^{\prime
}=B(\overline{y},\min\{\delta,\mu\}),V=B(\overline{z},\delta),W^{\prime
}=B(\overline{w},\delta).$ Take $z\in V$ and $(y,w)\in\operatorname*{Gr}%
G_{z}\cap\lbrack U^{\prime}\times W^{\prime}]$. Then $w\in(R_{y}+T)(z)\cap
B(\overline{w},\delta)$ and $(z,y)\in B(\overline{z},\delta)\times
B(\overline{y},\delta),$ whence there exist $w_{1}\in R(y)\cap B({\overline
{w}}_{1}{,}\varepsilon)$ and $w_{2}\in T(z)\cap B(\overline{w}_{2}%
,\varepsilon)$ such that $w=w_{1}+w_{2}.$ Consequently, $(y,w_{1}%
)\in\operatorname*{Gr}R\cap\lbrack U\times W_{1}],$ whence%
\[
B(w_{1},C\theta)\subset R(B(y,\theta)),
\]
i.e.%
\[
B(w_{1},C\theta)+w_{2}\subset R(B(y,\theta))+w_{2}.
\]
The shift-invariance of the distance in $W$ ensures
\[
B(w,C\theta)\subset R(B(y,\theta))+T(z)=G_{z}(B(y,\theta)).
\]
Therefore, we infer that $G$ is open at linear rate with respect to $y$
uniformly in $z$ around $((\overline{y},{\overline{z}),}\overline{w})$ with
constant $C>0.$

We have to prove now that $G$ has the Aubin property with respect to $z$
uniformly in $y$ around $((\overline{y},{\overline{z}),}\overline{w})$ with
constant $D>0.$ We know that there exists $V\in\mathcal{V}(\overline{z})$ and
$W_{2}\in\mathcal{V}(\overline{w}_{2})$ such that for all $z_{1},z_{2}\in V,$%
\[
e(T(z_{1})\cap W_{2},T(z_{2}))\leq Dd(z_{1},z_{2}).
\]
Take $\varepsilon>0$ such that $B(\overline{w}_{2},\varepsilon)\subset W_{2}.$
Take now $U=B(\overline{y},\delta),V^{\prime}=V\cap B(\overline{z},\delta)$
and $W^{\prime}=B(\overline{w},\delta),$ where, as above, $\delta$ is the
positive value given by Definition \ref{sum-sta_param} for the prescribed
$\varepsilon>0.$ Take $z_{1},z_{2}\in V^{\prime}$ and $w\in(R(y)+T(z_{1}))\cap
W^{\prime}.$ Then $w\in(R_{y}+T)(z_{1})\cap B(\overline{w},\delta)$ and
$(z_{1},y)\in B(\overline{z},\delta)\times B(\overline{y},\delta).$ Hence,
there exist $w_{1}\in R(y)\cap B({\overline{w}}_{1}{,}\varepsilon)$ and
$w_{2}\in T(z_{1})\cap B(\overline{w}_{2},\varepsilon)$ such that
$w=w_{1}+w_{2}.$ We deduce that $w_{2}\in T(z_{1})\cap W_{2},$ whence
$d(w_{2},T(z_{2}))\leq Dd(z_{1},z_{2}).$ Then
\begin{align*}
d(w,G_{y}(z_{2}))  &  =d(w_{1}+w_{2},R(y)+T(z_{2}))\\
&  \leq d(w_{1}+w_{2},w_{1}+T(z_{2}))=d(w_{2},T(z_{2}))\leq Dd(z_{1},z_{2}).
\end{align*}
Since $w\in(R_{y}+T)(z_{1})\cap B(\overline{w},\delta)$ was arbitrarily
chosen,%
\[
e(G_{y}(z_{1})\cap W^{\prime},G_{y}(z_{2}))\leq Dd(z_{1},z_{2})
\]
and the thesis is proved. The final conclusion is a direct application of
Theorem \ref{op_comp}.\hfill$\square$

\bigskip

In order to treat a similar situation in a different fashion we intend to
apply Theorem \ref{main_const} for $R\circ F_{1}$ and $-T\circ F_{2}.$ Since
the metric regularity of $R\circ F_{1}$\ take place provided that $R,$ $F_{1}$
share the same property (see the proof of Theorem \ref{t_sep_2}), we look now
at the Aubin property of $T\circ F_{2}$ and we present a definition.

\begin{df}
\label{lscomp}Let $F:X\rightrightarrows Y,$ $G:Y\rightrightarrows Z$ be
multifunctions and $(\overline{x},{\overline{y},}\overline{z})\in X\times
Y\times Z$ such that $\overline{y}\in F(\overline{x}),$ $\overline{z}\in
G(\overline{y}).$ We say that the pair of multifunctions $F,G$ is locally
composition-stable around $(\overline{x},\overline{y},\overline{z})$ if for
every $\varepsilon>0$ there exists $\delta>0$ such that, for every $x\in
B(\overline{x},\delta)$ and every $z\in(G\circ F)(x)\cap B(\overline{z}%
,\delta),$ there exists $y\in F(x)\cap B({\overline{y},}\varepsilon)$ such
that $z\in G(y).$
\end{df}

We provide next the main reason for introducing this notion.

\begin{lm}
\label{Ap_comp}Let $F:X\rightrightarrows Y,$ $G:Y\rightrightarrows Z$ be
multifunctions and $(\overline{x},{\overline{y},}\overline{z})\in X\times
Y\times Z$ such that $\overline{y}\in F(\overline{x}),$ $\overline{z}\in
G(\overline{y}).$ If $F$ and $G$ have the Aubin property around $(\overline
{x},{\overline{y}})\ $and $(\overline{y},\overline{z}),$ respectively, and
$F,G$ are locally composition-stable around $(\overline{x},{\overline{y}%
,}\overline{z}),$ then the multifunction $G\circ F$ has the Aubin property
around $(\overline{x},\overline{z}),$ and%
\[
\operatorname*{lip}(G\circ F)(\overline{x},\overline{z})\leq
\operatorname*{lip}F(\overline{x},{\overline{y}})\cdot\operatorname*{lip}%
G(\overline{y},\overline{z}).
\]

\end{lm}

\noindent\textbf{Proof.} According to the assumptions made on $F$ and $G,$ one
can find $\alpha>0,$ $l_{F}>0,l_{G}>0$ such that, for every $x,u\in
B(\overline{x},\alpha),$ and every $y,v\in B({\overline{y},}\alpha),$ one has%
\begin{align}
e(F(x)\cap B(\overline{y},\alpha),F(u))  &  \leq l_{F}d(x,u)\text{,}%
\label{AF}\\
e(G(y)\cap B(\overline{z},\alpha),G(v))  &  \leq l_{G}d(y,v). \label{AG}%
\end{align}

Applying the local stability of $F,G$ for $\varepsilon:=2^{-1}\alpha,$ one can
find $\delta\in(0,\min\{\alpha,(8l_{F})^{-1}\alpha\})$ such that the property
from Definition \ref{lscomp} is satisfied. Take now arbitrary $x,u\in
B(\overline{x},\delta)$ and $z\in(G\circ F)(x)\cap B(\overline{z},\delta).$
Then there exists $y\in F(x)\cap B({\overline{y},2}^{-1}\alpha)$ such that
$z\in G(y).$ So, according to (\ref{AF}), one has that%
\[
d(y,F(u))\leq l_{F}d(x,u).
\]

\noindent Hence, for every $\theta\in(0,4^{-1}\alpha),$ there exists $v\in
F(u)$ such that
\[
d(y,v)\leq l_{F}d(x,u)+\theta<\frac{\alpha}{2}.
\]
In conclusion, $y,v\in B({\overline{y},}\alpha),$ so one can apply (\ref{AG})
to get that%
\[
d(z,(G\circ F)(u))\leq d(z,G(v))\leq l_{G}d(y,v)\leq l_{G}l_{F}d(x,u)+\theta
l_{G}.
\]

But, because of the arbitrariness of $z\ $from $(G\circ F)(x)\cap
B(\overline{z},\delta),$ one gets that for every $x,u\in B(\overline{x}%
,\delta),$%
\[
e((G\circ F)(x)\cap B(\overline{z},\delta),(G\circ F)(u))\leq l_{G}%
l_{F}d(x,u)+\theta l_{G}.
\]

\noindent Making $\theta\rightarrow0,$ one gets the conclusion.\hfill$\square$

\bigskip

In the following, we provide an example of two multifunctions with the Aubin
property, but for which their composition fails to satisfy the same property.

\begin{examp}
Take $F:\mathbb{R\rightrightarrows R}$ given by
\[
F(x):=\left\{
\begin{array}
[c]{ll}%
\lbrack0,1], & \text{if }x\in\mathbb{R}\setminus\{1\}\\
\lbrack0,1]\cup\{2\}, & \text{if }x=1
\end{array}
\right.
\]
and $G:\mathbb{R\rightrightarrows R}$ given by
\[
G(y):=\left\{
\begin{array}
[c]{cc}%
\lbrack y,1), & \text{if }y<1\\
\{1\}, & \text{if }y=1\\
(1,y], & \text{if }y>1
\end{array}
\right.  .
\]

One can easily see that both $F$ and $G$ have the Aubin property around
$(\overline{x},\overline{y})=(1,1)$ and $(\overline{y},\overline{z})=(1,1),$
respectively. But the multifunction $G\circ F:\mathbb{R}\rightrightarrows
\mathbb{R}$, given by%
\[
(G\circ F)(x)=\left\{
\begin{array}
[c]{ll}%
\lbrack0,1], & \text{if }x\in\mathbb{R}\setminus\{1\}\\
\lbrack0,2], & \text{if }x=1
\end{array}
\right.  ,
\]

\noindent does not satisfy the same property around $(\overline{x}%
,\overline{z})=(1,1).$ Let us prove this last assertion: suppose by
contradiction that there exists $L>0$ and $\alpha\in(0,\min\{1,L^{-1}\})$ such
that for any $x,u\in\lbrack1-\alpha,1+\alpha]$%
\begin{equation}
(G\circ F)(x)\cap\lbrack1-\alpha,1+\alpha]\subset(G\circ F)(u)+L\left\vert
x-u\right\vert [-1,1]. \label{rpL}%
\end{equation}
Consider now $x:=1$ and $u:=1+\alpha^{2}$ such that $x,u\in\lbrack
1-\alpha,1+\alpha].$ Clearly,%
\[
1+\alpha\in(G\circ F)(x)\cap\lbrack1-\alpha,1+\alpha].
\]
Following (\ref{rpL}), we should have:%
\[
1+\alpha\in(G\circ F)(1+\alpha^{2})+L\alpha^{2}[-1,1]
\]
and, in particular,
\[
1+\alpha\in\lbrack-L\alpha^{2},1+L\alpha^{2}].
\]
But this requires that%
\[
\alpha\leq L\alpha^{2},
\]
which contradicts the choice of $\alpha$. The contradiction shows that we
cannot have the Aubin property of $G\circ F.$

Now, taking into account Lemma \ref{Ap_comp}, the pair $F,G$ cannot be locally
stable under composition around $(1,1,1).$ Indeed, pick $\varepsilon
\in(0,2^{-1}).$ Then for every $\delta>0,$ choose $n\in\mathbb{N}$ such that
$n>\max\{\delta,1\}.$ Taking now $x_{\delta}:=1\in(1-\delta,1+\delta)$ and
$z_{\delta}:=1+n^{-1}\delta\in(G\circ F)(x_{\delta})\cap(1-\delta,1+\delta),$
one can easily see that, for every $y\in F(x_{\delta})\cap(1-\varepsilon
,1+\varepsilon)=(1-\varepsilon,1],$ $z_{\delta}\not \in G(y).$
\end{examp}

\bigskip

\begin{rmk}
Let us observe that, if one takes in Definition \ref{lscomp}
$F:X\rightrightarrows Y\times Y,$ $F:=(F_{1},F_{2}),$ where $F_{1}%
:X\rightrightarrows Y,$ $F_{2}:X\rightrightarrows Y$ are two multifunctions,
$G:=g,$ where $g:Y\times Y\rightarrow Y$ is given by $g(y,z):=y+z,$ and
$(\overline{x},{\overline{y},}\overline{z})\in X\times Y\times Y$ such that
${\overline{y}\in F}_{1}{(\overline{x}),}$ ${\overline{z}\in F}_{2}%
{(\overline{x}),}$ then the local stability under composition of the pair
$F,G$ around $(\overline{x},({\overline{y},}\overline{z}),{\overline{y}%
+}\overline{z})$ is just the local sum-stability of $(F_{1},F_{2})$ around
$(\overline{x},{\overline{y},}\overline{z}).$ Also, in view of Lemma
\ref{Ap_comp}, one gets that the sum of two multifunctions with the Aubin
property around corresponding points has the Aubin property provided that the
two multifunctions are locally sum-stable. For more details in this direction,
see \cite{DurStr4}, Definition 4.2 and the subsequent examples and results.
\end{rmk}

\bigskip

Putting all these facts together, we are now in position to formulate another
theorem for the situation of separate variables.

\begin{thm}
\label{t_sep_2}Let $X,Y,Z$ be metric spaces, $W$ be a metric space with linear
structure such that the metric is shift-invariant. Take $F_{1}%
:X\rightrightarrows Y,$ $F_{2}:X\rightrightarrows Z,$ $R:Y\rightrightarrows
W,$ $T:Z\rightrightarrows W$ and $(\overline{x},\overline{y},\overline
{z},\overline{w}_{1},\overline{w}_{2})\in X\times Y\times Z\times W\times W$
such that $\overline{y}\in F_{1}(\overline{x}),$ $\overline{w}_{1}\in
R(\overline{y}),$ $\overline{z}\in F_{2}(\overline{x}),$ $\overline{w}_{2}\in
T(\overline{z})$. Suppose that:

(i) $\operatorname*{Gr}(R\circ F_{1})$ and $\operatorname*{Gr}(T\circ F_{2})$
are locally complete around $(\overline{x},\overline{w}_{1})$ and
$(\overline{x},\overline{w}_{2}),$ respectively.

(ii) $F_{1}$ is open at linear rate $L>0$ around $(\overline{x},\overline
{y});$

(iii) $R$ is open at linear rate $C>0$ around $(\overline{y},\overline{w}%
_{1});$

(iv) $F_{2}$ has the Aubin property around $(\overline{x},\overline{z})$ with
constant $M>0;$

(v) $T$ has the Aubin property around $(\overline{z},\overline{w}_{2})$ with
constant $D>0;$

(vi) $F_{2},T$ are locally composition-stable at $(\overline{x},\overline
{z},\overline{w}_{2});$

(vii) $LC-MD>0.$

Then there exists $\varepsilon>0$ such that for every $\rho\in(0,\varepsilon)$
such that%
\[
B(\overline{w}_{1}+\overline{w}_{2},(LC-MD)\rho)\subset(R\circ F_{1}+T\circ
F_{2})(B(\overline{x},\rho)).
\]

Moreover, there exists $\varepsilon^{\prime}>0$ such that for every $\rho
\in(0,\varepsilon^{\prime})$ and every $(x,y,z,w_{1},w_{2})\in B(\overline
{x},\varepsilon^{\prime})\times B(\overline{y},\varepsilon^{\prime})\times
B(\overline{z},\varepsilon^{\prime})\times B(\overline{w}_{1},\varepsilon
^{\prime})\times B(\overline{w}_{2},\varepsilon^{\prime})$ such that $y\in
F_{1}(x),$ $w_{1}\in R(y),$ $z\in F_{2}(x),$ $w_{2}\in T(z)$%
\[
B(w_{1}+w_{2},(LC-MD)\rho)\subset(R\circ F_{1}+T\circ F_{2})(B(x,\rho)).
\]

\end{thm}

\noindent\textbf{Proof.} Remark that, using (ii) and (iii), the multifunction
$R\circ F_{1}$ is $LC-$open around $(\overline{x},\overline{w}_{1}).$ Also,
from (iv), (v) and (vi), using Lemma \ref{Ap_comp}, one gets that $T\circ
F_{2}$ has the Aubin property around $(\overline{x},\overline{w}_{2})$ with
constant $MD.$ Next, consider $G$ as in Remark \ref{rmk_sum} and apply Theorem
\ref{main_const} for $R\circ F_{1}$ and $-T\circ F_{2}.$\hfill$\square$

\bigskip

Remark that the main differences between Theorems \ref{t_sep_1} and
\ref{t_sep_2} are, on one side, those refering to the completeness and the
closedness of the graphs, and, on the other side, those concerning the local stability.

\subsection{Parametric variational systems}

In the sequel, we shall need a result previously given in \cite{DurStr4},
which makes the link between a parametric multifunction and the associated
solution map, providing also interesting metric evaluations and relations
between the regularity moduli of involved set-valued mappings. To this aim,
consider a multifunction $H:X\times P\rightrightarrows W,$ where $X,P$ are
metric spaces, and $W$ is a normed vector space and define the implicit
solution map $S:P\rightrightarrows X$ by%
\[
S(p):=\{x\in X\mid0\in H(x,p)\}.
\]

The next implicit multifunction theorem will play an important role since it
will provide both metric regularity and Aubin property for $S$. The full
version of this result is done in \cite[Theorem 3.6]{DurStr4}.

\begin{thm}
\label{impl} Let $X,P$ be metric spaces, $Y$ be a normed vector space,
$H:X\times P\rightrightarrows W$ be a set-valued map and $(\overline
{x},\overline{p},0)\in\operatorname{Gr}H$.

(i) If $H$ is open at linear rate $c>0$ with respect to $x$ uniformly in $p$
around $(\overline{x},\overline{p},0)$, then there exist $\alpha,\beta
,\gamma>0$ such that, for every $(x,p)\in B(\overline{x},\alpha)\times
B(\overline{p},\beta),$%
\begin{equation}
d(x,S(p))\leq c^{-1}d(0,H(x,p)\cap B(0,\gamma)). \label{xSpV}%
\end{equation}

Suppose, in addition to (\ref{xSpV}), that $H$ has the Aubin property with
respect to $p$ uniformly in $x$ around $(\overline{x},\overline{p},0).$ Then
$S$ has the Aubin property around $(\overline{p},\overline{x})$ and%
\begin{equation}
\operatorname*{lip}S(\overline{p},\overline{x})\leq c^{-1}\widehat
{\operatorname*{lip}}_{p}H((\overline{x},\overline{p}),0). \label{lipS}%
\end{equation}

(ii) If $H$ is open at linear rate $c>0$ with respect to $p$ uniformly in $x$
around $(\overline{x},\overline{p},0)$, then there exist $\alpha,\beta
,\gamma>0$ such that, for every $(x,p)\in B(\overline{x},\alpha)\times
B(\overline{p},\beta),$%
\begin{equation}
d(p,S^{-1}(x))\leq c^{-1}d(0,H(x,p)\cap B(0,\gamma)). \label{pSxV}%
\end{equation}

Suppose, in addition to (\ref{pSxV}), that $H$ has the Aubin property with
respect to $x$ uniformly in $p$ around $(\overline{x},\overline{p},0).$ Then
$S$ is metrically regular around $(\overline{p},\overline{x})$ and%
\begin{equation}
\operatorname*{reg}S(\overline{p},\overline{x})\leq c^{-1}\widehat
{\operatorname*{lip}}_{x}H((\overline{x},\overline{p}),0). \label{regS}%
\end{equation}

\end{thm}

\bigskip

Take $H(x,p):=G(F_{1}(x),F_{2}(x,p)),$ with $F_{1}:X\rightrightarrows Y,$
$F_{2}:X\times P\rightrightarrows Z,$ $G:Y\times Z\rightrightarrows W.$ Then
the openness result in Theorem \ref{op_comp} and the previous implicit
multifunction theorem come into play, to ensure results concerning the
well-posedness of the solution mapping associated to the next parametric
variational system%
\begin{equation}
0\in G(F_{1}(x),F_{2}(x,p)). \label{PVS}%
\end{equation}

\bigskip

Also, for two multifunctions $F_{1}:X\rightrightarrows Y,$ $F_{2}:X\times
P\rightrightarrows Z,$ we consider (as in Theorem \ref{t_sep_1} (vi)) the
multifunction $(F_{1},F_{2}):X\times P\rightrightarrows Y\times Z$ given by%
\[
(F_{1},F_{2})(x,p):=F_{1}(x)\times F_{2}(x,p).
\]

We are now in position to formulate our results concerning the metric
regularity and the Aubin property of the solution mapping associated to
(\ref{PVS}).

Practically, we follow the same way as in \cite[Theorems 4.12, 4.13]{DurStr4},
this time on more general setting and parametric systems. This approach was
recently brought into attention by the works of Dontchev and Rockafellar
\cite{DontRock2009a} and Arag\'{o}n Artacho and Mordukhovich
\cite{ArtMord2009}, \cite{ArtMord2010}. Finally, let us remark that the
estimations we obtain here cover those in the quoted papers.

\begin{thm}
\label{mreg_sol}Let $X,P,Y,Z$ be metric spaces, $W$ be a normed vector space,
$F_{1}:X\rightrightarrows Y,$ $F_{2}:X\times P\rightrightarrows Z,$ $G:Y\times
Z\rightrightarrows W$ be set-valued maps and $(\overline{x},\overline
{p},\overline{y},\overline{z})\in X\times P\times Y\times Z$ such that
$\overline{y}\in F_{1}(\overline{x}),$ $\overline{z}\in F_{2}(\overline
{x},\overline{p})$ and $0\in G(\overline{y},\overline{z})$. Suppose that the
following assumptions are satisfied:

(i) $(F_{1},F_{2}),G$ are locally composition-stable around $((\overline
{x},\overline{p}),({\overline{y},\overline{z}),0});$

(ii) $F_{1}$ has the Aubin property around $(\overline{x},\overline{y});$

(iii) $F_{2}$ has the Aubin property with respect to $x$ uniformly in $p$
around $((\overline{x},\overline{p}),\overline{z});$

(iv) $F_{2}$ is metrically regular with respect to $p$ uniformly in $x$ around
$((\overline{x},\overline{p}),\overline{z});$

(v) $G$ is metrically regular with respect to $z$ uniformly in $y$ around
$((\overline{y},\overline{z}),0);$

(vi) $G$ has the Aubin property around $((\overline{y},\overline{z}),0).$

Then $S$ is metrically regular around $(\overline{p},\overline{x}).$ Moreover,
the next relation holds%
\begin{equation}
\operatorname*{reg}S(\overline{p},\overline{x})\leq\widehat
{\operatorname*{reg}}_{p}F_{2}((\overline{x},\overline{p}),\overline{y}%
)\cdot\widehat{\operatorname*{reg}}_{z}G((\overline{y},\overline{z}%
),0)\cdot\max\{\operatorname*{lip}F_{1}(\overline{x},\overline{z}%
),\widehat{\operatorname*{lip}}_{x}F_{2}((\overline{x},\overline{p}%
),\overline{y})\}\cdot\operatorname*{lip}G((\overline{y},\overline{z}),0).
\label{rS}%
\end{equation}

\end{thm}

\noindent\textbf{Proof. }Consider the multifunction $H:X\times
P\rightrightarrows W$ given by%
\begin{equation}
H(x,p):=(G\circ(F_{1},F_{2}))(x,p). \label{H}%
\end{equation}

Using (ii) and (iii), one can easily prove that $(F_{1},F_{2})$ has the Aubin
property with respect to $x$ uniformly in $p$ around $((\overline{x}%
,\overline{p}),{(\overline{y},\overline{z})})$ with modulus $K:=\max
\{\operatorname*{lip}F_{1}(\overline{x},\overline{z}),\widehat
{\operatorname*{lip}}_{x}F_{2}((\overline{x},\overline{p}),\overline{y})\}.$
In view of (i), (vi) and Lemma \ref{Ap_comp}, we know that $H$ has the Aubin
property with respect to $x$ uniformly in $p$ around $((\overline{x}%
,\overline{p}),{0})$ and the relation%
\[
\widehat{\operatorname*{lip}}_{x}H((\overline{x},\overline{p}),0)\leq
K\cdot\operatorname*{lip}G((\overline{y},\overline{z}),0)
\]
holds.

Using now Theorem \ref{link_around}, (iv) is equivalent to the fact that
$F_{2}$ is open at linear rate with respect to $p$ uniformly in $x$ around
$((\overline{x},\overline{p}),\overline{z})$ and $\widehat{\operatorname*{lop}%
}_{p}F_{2}((\overline{x},\overline{p}),\overline{y})=(\widehat
{\operatorname*{reg}}_{p}F_{2}((\overline{x},\overline{p}),\overline{y}%
))^{-1}.$ Similarly, $G$ is open at linear rate with respect to $z$ uniformly
in $y$ around $((\overline{y},\overline{z}),0)$ and $\widehat
{\operatorname*{lop}}_{z}G((\overline{y},\overline{z}),0)=(\widehat
{\operatorname*{reg}}_{z}G((\overline{y},\overline{z}),0))^{-1}.$

Consequently, there exist $\varepsilon,L,C>0$ such that, for every
$(x,p,y,z,w)\in B(\overline{x},\varepsilon)\times B(\overline{p}%
,\varepsilon)\times B(\overline{y},\varepsilon)\times B(\overline
{z},\varepsilon)\times B(\overline{w},\varepsilon)$ such that $(p,z)\in
\operatorname*{Gr}(F_{2})_{x}$ and $(z,w)\in\operatorname*{Gr}G_{y},$ and
every $\rho\in(0,\varepsilon),$%
\begin{align}
B(z,\rho)  &  \subset(F_{2})_{x}(B(p,L^{-1}\rho)),\label{opF1}\\
B(w,C\rho)  &  \subset G_{y}(B(z,\rho)). \label{op_G}%
\end{align}

Using now the local stability from (i), there exists $\delta\in(0,\varepsilon
)$ such that, for every $(x,p)\in B(\overline{x},\delta)\times B(\overline
{p},\delta)$ and every $w\in(G\circ(F_{1},F_{2}))(x,p)\cap B(0,\delta),$ there
exists $(y,z)\in\lbrack F_{1}(x)\times F_{2}(x,p)]\cap\lbrack B({\overline
{y},}\varepsilon)\times B(\overline{z},\varepsilon)]$ such that $w\in G(y,z).$

Take now arbitrary $x\in B(\overline{x},\delta),$ $(p,w)\in\operatorname*{Gr}%
H_{x}\cap\lbrack B(\overline{p},\delta)\times B(0,\delta)]$ and $\rho
\in(0,\varepsilon).$ Then $w\in H(x,p)\cap B(0,\delta),$ so there exists
$(y,z)$ as above such that $(p,z)\in\operatorname*{Gr}(F_{2})_{x}\cap\lbrack
B(\overline{p},\varepsilon)\times B(\overline{z},\varepsilon)].$ Also,
$(z,w)\in\operatorname*{Gr}G_{y}\cap\lbrack B(\overline{z},\varepsilon)\times
B(\overline{w},\varepsilon)].$ Hence, using (\ref{opF1}) and (\ref{op_G}),
\[
B(w,C\rho)\subset G(y,B(z,\rho))\subset G(y,F_{2}(x,B(p,L^{-1}\rho)))\subset
H(x,B(p,L^{-1}\rho)).
\]

In conclusion, $H$ is open at linear rate with respect to $p$ uniformly in
$x,$ and $\widehat{\operatorname*{lop}}_{p}H((\overline{x},\overline
{p}),0)\leq\widehat{\operatorname*{lop}}_{p}F_{2}((\overline{x},\overline
{p}),\overline{y})\cdot\widehat{\operatorname*{lop}}_{z}G((\overline
{y},\overline{z}),0)=(\widehat{\operatorname*{reg}}_{p}F_{2}((\overline
{x},\overline{p}),\overline{y}))^{-1}\cdot(\widehat{\operatorname*{reg}}%
_{z}G((\overline{y},\overline{z}),0))^{-1}.$

Now the result follows from Theorem \ref{impl} (ii).$\hfill\square$

\bigskip

Next, we present a more involved result, which makes use of the Theorem
\ref{op_comp}.

\bigskip

\begin{thm}
\label{Lip_sol}Let $X,P,Y,Z$ be metric spaces, $W$ be a normed vector space,
$F_{1}:X\rightrightarrows Y,$ $F_{2}:X\times P\rightrightarrows Z,$ $G:Y\times
Z\rightrightarrows W$ be set-valued maps and $(\overline{x},\overline
{p},\overline{y},\overline{z})\in X\times P\times Y\times Z$ such that
$\overline{y}\in F_{1}(\overline{x}),$ $\overline{z}\in F_{2}(\overline
{x},\overline{p})$ and $0\in G(\overline{y},\overline{z})$. Suppose that the
following assumptions are satisfied:

(i) $(F_{1},F_{2}),G$ is locally composition-stable around $((\overline
{x},\overline{p}),({\overline{y},}\overline{z}),0);$

(ii) $\operatorname*{Gr}F_{1}$ is complete, $\operatorname*{Gr}(F_{2})_{p}$ is
complete for every $p$ in a neighborhood of $\overline{p}$, and
$\operatorname*{Gr}G$ is closed;

(iii) $F_{1}$ is open at linear rate around $(\overline{x},\overline{y});$

(iv) $F_{2}$ has the Aubin property around $((\overline{x},\overline
{p}),\overline{z});$

(v) $G$ is open at linear rate with respect to $y$ uniformly in $z$ around
$((\overline{y},\overline{z}),0);$

(vi) $G$ has the Aubin property with respect to $z$ uniformly in $y$ around
$((\overline{y},\overline{z}),0);$

(vii) $\widehat{\operatorname*{lip}}_{x}F_{2}((\overline{x},\overline
{p}),\overline{y})\cdot\widehat{\operatorname*{lip}}_{z}G((\overline
{y},\overline{z}),0)<\widehat{\operatorname*{lop}}_{y}G((\overline
{y},\overline{z}),0)\cdot\operatorname*{lop}F_{1}(\overline{x},\overline{z}).$

Then $S$ has the Aubin property around $(\overline{p},\overline{x})$.
Moreover, the next relation is satisfied%
\begin{equation}
\operatorname*{lip}S(\overline{p},\overline{x})\leq\frac{\widehat
{\operatorname*{lip}}_{p}F_{2}((\overline{x},\overline{p}),\overline{y}%
)\cdot\widehat{\operatorname*{lip}}_{z}G(({\overline{y},}\overline{z}%
),0)}{\widehat{\operatorname*{lop}}_{y}G((\overline{y},\overline{z}%
),0)\cdot\operatorname*{lop}F_{1}(\overline{x},\overline{z})-\widehat
{\operatorname*{lip}}_{x}F_{2}((\overline{x},\overline{p}),\overline{y}%
)\cdot\widehat{\operatorname*{lip}}_{z}G((\overline{y},\overline{z}),0)}.
\label{lS}%
\end{equation}

\end{thm}

\noindent\textbf{Proof. }Take $L>\operatorname*{lop}F_{1}(\overline
{x},\overline{y}),$ $C>\widehat{\operatorname*{lop}}_{y}G((\overline
{y},\overline{z}),0),$ $M>\widehat{\operatorname*{lip}}_{x}F_{2}((\overline
{x},\overline{p}),\overline{y})$ and $D>\widehat{\operatorname*{lip}}%
_{z}G((\overline{y},\overline{z}),0)$ such that $LC-MD>0.$

Now, we intend to prove that there exist $\tau,t,\gamma>0$ such that, for
every $(x,p)\in B(\overline{x},\tau)\times B(\overline{p},t),$%
\begin{equation}
d(x,S(p))\leq(LC-MD)^{-1}d(0,H(x,p)\cap B(0,\gamma)). \label{mrg}%
\end{equation}

Using assumptions (ii)-(vi), one can find $\alpha>0$ such that:

\begin{enumerate}
\item $\operatorname{Gr}F_{1}\cap\lbrack D(\overline{x},\alpha)\times
D({\overline{y}},\alpha)]$ is complete; for every $p\in B(\overline{p}%
,\alpha),$ $\operatorname{Gr}(F_{2})_{p}\cap\lbrack D(\overline{x}%
,\alpha)\times D(\overline{z},\alpha)]$ is complete; $\operatorname{Gr}%
G\cap\lbrack D(\overline{y},\alpha)\times D(\overline{z},\alpha)\times
D(0,\alpha)]$ is closed.

\item for every $(x,y)\in B(\overline{x},\alpha)\times B({\overline{y}}%
,\alpha),$
\begin{equation}
d(x,F_{1}^{-1}(y))\leq\frac{1}{L}d(y,F_{1}(x)) \label{lopF1}%
\end{equation}

\item for every $p\in B(\overline{p},\alpha)$ and every $x,x^{\prime}\in
B({\overline{x},}\alpha),$%
\begin{equation}
e(F_{2}(x,p)\cap B(\overline{z},\alpha),F_{2}(x^{\prime},p))\leq
Md(x,x^{\prime}). \label{lip_F2}%
\end{equation}

\item for every $(z,w),(z^{\prime},w^{\prime})\in B({\overline{z},}%
\alpha)\times B(0,\alpha)$,%
\begin{equation}
e(\Gamma(z,w)\cap D(\overline{y},\alpha),\Gamma(z^{\prime},w^{\prime}%
))\leq\frac{1}{C}(Dd(z,z^{\prime})+d(w,w^{\prime})\mathbb{)}. \label{rel_Gam}%
\end{equation}

\item for every $y\in B(\overline{y},\alpha)$ and every $z,z^{\prime}\in
B({\overline{z},}\alpha),$%
\begin{equation}
e(G(y,z)\cap B(0,\alpha),G(y,z^{\prime}))\leq Dd(z,z^{\prime}). \label{lip_G}%
\end{equation}

\end{enumerate}

Choose now $\varepsilon>0$ such that (\ref{eps}) are satisfied with
$2^{-1}\alpha$ instead of $\alpha.$ Finally, apply the property from (i) for
$2^{-1}\alpha$ instead of $\varepsilon$ and find $\delta\in(0,2^{-1}\alpha)$
such that the assertion from Definition \ref{lscomp} is true.

Take now $\rho\in(0,\min\{(LC-MD)^{-1}\delta,\varepsilon\}),\ $define
$\gamma:=(LC-MD)\rho$ and fix $(x,p)\in B(\overline{x},\delta)\times
B(\overline{p},\delta).$

If $H(x,p)\cap B(0,\gamma)=\emptyset$ or $0\in H(x,p)\cap B(0,\gamma),$ then
(\ref{mrg}) trivially holds. Suppose next that $0\not \in H(x,p)\cap
B(0,\gamma).$ Then, for every $\theta>0,$ one can find $w_{\theta}\in
H(x,p)\cap B(0,\gamma)$ such that%
\begin{equation}
\left\Vert w_{\theta}\right\Vert <d(0,H(x,p)\cap B(0,\gamma))+\theta.
\label{ineg}%
\end{equation}

Because $d(0,H(x,p)\cap B(0,\gamma))<(LC-MD)\rho,$ for sufficiently small
$\theta,$ $d(0,H(x,p)\cap B(0,\gamma))+\theta<(LC-MD)\rho.$ Hence, it follows
from (\ref{ineg}) that%
\begin{equation}
0\in B(w_{\theta},d(0,H(x,p)\cap B(0,\gamma))+\theta)\subset B(w_{\theta
},(LC-MD)\rho)\subset B(w_{\theta},\delta), \label{0in}%
\end{equation}

\noindent so $w_{\theta}\in H(x,p)\cap B(0,\delta).$ Because we have also that
$(x,p)\in B(\overline{x},\delta)\times B(\overline{p},\delta),$ one can apply
(i) to find $y_{\theta}\in F_{1}(x)\cap B({\overline{y},2}^{-1}\alpha)$ and
$z_{\theta}\in F_{2}(x,p)\cap B(\overline{z}{,2}^{-1}\alpha)$ such that
$w_{\theta}\in G(y_{\theta},z_{\theta}).$ Consequently, $B(y_{\theta}%
,2^{-1}\alpha)\subset B(\overline{y},\alpha)$ and $B(z_{\theta},2^{-1}%
\alpha)\subset B(\overline{z},\alpha).$

Observe now that the relations (\ref{lopF1})-(\ref{lip_G}) are satisfied for
$x,p,y_{\theta},z_{\theta},w_{\theta}$ instead of $\overline{x},\overline
{p},\overline{y},\overline{z},0$ and $2^{-1}\alpha$ instead of $\alpha,$
because every ball centered in these points with radius $2^{-1}\alpha$ is
contained in the initial one with radius $\alpha.$

Now, because $\varepsilon$ was chosen such that (\ref{eps}) is satisfied for
$2^{-1}\alpha$ instead of $\alpha,$ one can use Theorem \ref{op_comp} for
$F_{1},(F_{2})_{p},G,$ the reference points $y_{\theta}\in F_{1}(x),$
$z_{\theta}\in(F_{2})_{p}(x),$ $w_{\theta}\in G(y_{\theta},z_{\theta}),$ and
for $\rho_{0}:=(LC-MD)^{-1}\cdot(d(0,H(x,p)\cap B(0,\gamma))+\theta
)<\rho<\varepsilon$ to get that%
\[
B(w_{\theta},d(0,H(x,p)\cap B(0,\gamma))+\theta)\subset G\circ(F_{1}%
,(F_{2})_{p})(B(x,\rho_{0})).
\]

Using also (\ref{0in}), we know that $0\in G\circ(F_{1},(F_{2})_{p}%
)(B(x,\rho_{0})),$ so there exists $\widetilde{x}\in B(x,\rho_{0})$ such that
$0\in G(F_{1}(\widetilde{x}),F_{2}(\widetilde{x},p))$ or, equivalently,
$\widetilde{x}\in S(p).$ Consequently,%
\[
d(x,S(p))\leq d(x,\widetilde{x})<\rho_{0}=(LC-MD)^{-1}\cdot(d(0,H(x,p)\cap
B(0,\gamma))+\theta).
\]

Making $\theta\rightarrow0,$ one gets (\ref{mrg}).

Now, for the final step of the proof, observe that, from the Aubin property of
$F_{2}$ with respect to $p$ uniformly in $x$ around $((\overline{x}%
,\overline{p}),\overline{z}),$ one can find $\beta,k>0$ such that for every
$x\in B(\overline{x},\beta),$ every $p_{1},p_{2}\in B(\overline{p},\beta),$
one has%
\begin{equation}
e(F_{2}(x,p_{1})\cap B(\overline{z},\beta),F_{2}(x,p_{2}))\leq kd(p_{1}%
,p_{2}). \label{AubF2}%
\end{equation}

Denote $\xi:=\min\{\alpha,\beta\}.$ One can use now (i) for $2^{-1}\xi$
instead of $\varepsilon$ to find $\delta^{\prime}\in(0,\min\{\xi,\delta
,6^{-1}k^{-1}\xi\})$ such that the assertion from Definition \ref{lscomp} is
true. Take now arbitrary $x\in B(\overline{x},\delta^{\prime}),$ $p_{1}%
,p_{2}\in B(\overline{p},\delta^{\prime})$ and $w\in H(x,p_{1})\cap
B(0,\delta^{\prime}).$ Then there exist $y\in F_{1}(x)\cap B(\overline
{y},2^{-1}\xi)$ and $z\in F_{2}(x,p_{1})\cap B(\overline{z},2^{-1}\xi)$ such
that $w\in G(y,z).$ Using now (\ref{AubF2}), one obtains that%
\[
d(z,F_{2}(x,p_{2}))\leq kd(p_{1},p_{2}).
\]

Hence, for every $\mu>0,$ there exists $z_{\mu}\in F_{2}(x,p_{2})$ such that
$d(z,z_{\mu})<kd(p_{1},p_{2})+\mu.$ Then%
\[
d(z_{\mu},\overline{z})\leq d(z_{\mu},z)+d(z,\overline{z})<kd(p_{1},p_{2}%
)+\mu+2^{-1}\xi<2k\delta^{\prime}+\mu+2^{-1}\xi<6^{-1}5\xi+\mu.
\]
But this means, for sufficiently small $\mu,$ that $z_{\mu}\in B(\overline
{z},\xi).$ Consequently, because of (\ref{lip_G}) and taking into account that
$w\in G(y,z)\cap B(0,\delta^{\prime}),$ we deduce that%
\[
d(w,G(y,z_{\mu}))\leq Dd(z,z_{\mu})<Dkd(p_{1},p_{2})+D\mu.
\]

Finally, because $y\in F_{1}(x)$ and $z_{\mu}\in F_{2}(x,p_{2}),$ one gets
that $d(w,H(x,p_{2}))\leq d(w,G(y,z_{\mu})),$ so%
\[
d(w,H(x,p_{2}))<Dkd(p_{1},p_{2})+D\mu.
\]

Making $\mu\rightarrow0$ and taking into account the arbitrariness of $w$ from
$H(x,p_{1})\cap B(0,\delta^{\prime}),$ we obtain that $H$ has the Aubin
property with respect to $p$ uniformly in $x$ around $((\overline{x}%
,\overline{p}),0).$ Moreover, $\widehat{\operatorname*{lip}}_{p}%
H((\overline{x},\overline{p}),0)\leq\widehat{\operatorname*{lip}}_{p}%
F_{2}((\overline{x},\overline{p}),\overline{y})\cdot\widehat
{\operatorname*{lip}}_{z}G(({\overline{y},}\overline{z}),0).$

Now the final conclusion follows from Theorem \ref{impl} (i).$\hfill\square$

\end{document}